\documentclass[reqno]{amsart}
\usepackage{amssymb}
\usepackage{amsmath}

\makeatletter
\@addtoreset{equation}{section}
\makeatother

\renewcommand\thefigure{\thesection.\@arabic\c@figure}
\renewcommand\thetable{\thesection.\@arabic\c@table}

\newtheorem{theorem}{Theorem}[section]
\newtheorem{lemma}[theorem]{Lemma}
\newtheorem{proposition}[theorem]{Proposition}

\newcommand{\mc}[1]{{\mathcal #1}}
\newcommand{\mf}[1]{{\mathfrak #1}}
\newcommand{\mb}[1]{{\mathbf #1}}
\newcommand{\bb}[1]{{\mathbb #1}}
\newcommand{\bs}[1]{{\boldsymbol #1}}

\begin{document}

\title{Superdiffusivity of two dimensional lattice gas models}

\author[C. Landim, J. A. Ram\'{\i}rez, H.-T. Yau]
{Claudio Landim$^\dagger$, Jos\'e  A. Ram\'{\i}rez$^\ddag$,
Horng-Tzer Yau$^\S$}

\thanks{$\dagger$ IMPA and CNRS Rouen, $\ddag$ Universidad de Costa Rica, $\S$ Stanford
  University on leave from Courant Institute, New York University}

\address{\noindent IMPA, Estrada Dona Castorina 110,
CEP 22460 Rio de Janeiro, Brasil and CNRS UMR 6085,
Universit\'e de Rouen, 76128 Mont Saint Aignan, France.
\newline
e-mail:  \rm \texttt{landim@impa.br}
}

\address{\noindent Escuela de Matem\'aticas, Universidad de Costa Rica, San Jos\'e 2060, Costa Rica.
\newline
e-mail:  \rm \texttt{jramirez@emate.ucr.ac.cr}
}

\address{\noindent Department of Mathematics, Stanford University, CA-94305,
USA.
\newline
e-mail: \rm \texttt{yau@math.stanford.edu}
}

\begin{abstract}It was proved \cite{EMYa, QY} that
stochastic lattice gas dynamics converge to the Navier-Stokes
equations in dimension $d=3$ in the incompressible limits. In
particular, the viscosity is finite. We proved that, on the other
hand,  the viscosity for a two dimensional lattice gas model
diverges faster than $\log \log t$. Our argument indicates that
the correct divergence rate is $(\log t)^{1/2}$. This problem is
closely related to the logarithmic correction of the time decay
rate for the velocity auto-correlation function of a tagged
particle.

\end{abstract}

\subjclass[2000]{primary 60K35}

\keywords{Hydrodynamic limit, second class particle, superdiffusivity}

\maketitle

\section{Introduction}
\label{sec0}

It is well-known that although the classical dynamics is time
reversible, the macroscopic behavior of the fluid, governed by the
Navier-Stokes equations, is time irreversible.  The measure on the
time irreversibility is characterized by the viscosity, also called
the bulk diffusivity of the system.  It can be represented as the
diffusion coefficient of a second class particle in the fluid.
Instead of a second class particle, one can study the typical behavior
of a tagged particle.  Once again, even though the underlying dynamics
is time reversible, the tagged particle is diffusive. The diffusion
coefficient in this case is called the self diffusion coefficient. The
bulk and self diffusion coefficients are two different quantities, but
they share similar qualitative behavior.

The Green-Kubo formulae represent the bulk or self diffusion
coefficients as time integrals over the current correlation function
or the velocity correlation function. In the fundamental work of Alder
and Wainwright \cite{AW}, it predicts that the time decay of the
velocity correlation function is of order $t^{-d/2}$, here $d$ is the
dimension of the system. Since the decay is $1/t$ in $d=2$ is not
integrable, it predicted that the self-diffusion coefficient of the
two dimensional fluid diverges. In the later work, Alder, Wainwright
and Gass \cite{AWG} proposed that the decay in two dimension is
actually $t^{-1} {\log t}^{-1/2}$. However, the logarithmic correction
cannot be seen from the their simulation.  This decay rate was also
obtained by Forster, Nelson and Stephen \cite{FNS} by a
renormalization group method.  The later simulation by van der Hoef
and Frenkel \cite{HF} confirmed that there is a discrepancy to the
pure algebraic decay $t^{-1}$, but is far from being able to determine
the precise logarithmic correction.

If we formally integrate the law $t^{-1} (\log t)^{-1/2}$, we obtain
from the Green-Kubo formula that the diffusion coefficient diverges as
$ (\log t)^{1/2}$ in dimension $d=2$. For higher dimension, the
diffusion coefficient is expected to be finite. This was proved
rigorously in various settings. For the stochastic lattice gas models
considered in \cite{EMYa}, the bulk diffusion coefficient is proved to
be finite for $d\ge 3$. One key ingredient of these lattice gas models
is the asymmetric simple exclusion process. For this process, both the
bulk and self diffusion coefficients were proved to be finite in
\cite{LY} and \cite{SVY} for dimension $d\ge 3$.  Van Beijeren, Kutner
and Spohn \cite{BKS} predicted via the mode-coupling theory that the
diffusivity diverges as $(\log t)^{2/3}$ in $d=2$ and $t^{1/3}$ for
$d=1$. The two dimensional case was proved in \cite{Y}; the one
dimensional case was only partly solved \cite{LQSY}. However, related
problems in the one dimensional case was solved by integrable method
\cite{J}.

In this paper, we shall proved that for the stochastic lattice model,
the bulk diffusion coefficient diverges in $d=2$. Following the method
of \cite{LQSY}, we derive a series of upper and lower bounds to the
diffusivity in terms of continuous fractions for operators.  If we
take the first lower bound to the diffusivity, we obtain the
divergence rate in Theorem~\ref{s0}. If one assumes that the
dispersion laws of these two sequences (upper and lower) of the
continuous fractions converge, the divergence law $(\log t)^{1/2}$ can
be obtained heuristically.  See the discussion at the end of this
paper.  Notice that the exponent $1/2$ is different from the $2/3$ in
the case of the asymmetric simple exclusion process.  Since the
lattice gas models are very complicated, we do not know if the
argument of \cite{Y} can be extended to this case.

\section{The model}
\label{sec1}

We recall the lattice gas models considered in \cite{EMYa} in
dimension $d$.  Denote by $\{e_j, j=1,\dots ,d\}$ the canonical basis
of $\bb R^d$ and let $\mc E = \{\pm e_1,\dots, \pm e_d\}$.  Let $
\mc V \subset \bb R^d $ be a finite set representing the possible
velocities.  On each site of the lattice at most one particle for each
velocity is allowed.  A configuration of particles on the lattice is
denoted by $\eta= \{\eta_x,~x\in \bb Z^d\}$ where
$\eta_x=\{\eta(x,v),~v\in{\mc V}\}$ and
$\eta(x,v)\in\{0,1\},~x\in\bb Z^d,~v\in{\mc V}$, is the number of
particles with velocity $v$ at $x$.  The set of particle
configurations is $X=\left(\{0,1\}^{\mc V}\right)^{\bb Z^d}$.

The dynamics consists of two parts: Asymmetric random walk with
exclusion among particles of the same velocity and binary collisions
between particles of different velocities.  We first describe the
random walk part of the dynamics.  Particles of velocity $v$ perform a
continuous time asymmetric random walk with simple exclusion.  A
particle at $x$ waits a random, exponentially distributed time then
chooses a nearby site $x+y$ according to a certain jump law and jumps
there as long as the site is not occupied by another particle of the
same velocity.  If there is a particle of the same velocity, the jump
is suppressed and the particle waits for a new exponential time.  All
particles are doing this simultaneously, and since time is continuous
ties do not occur. The jump law and waiting time are chosen so that
the rate of jumping from site $x$ to site $x+y$ is $p(y,v)$ which
should be finite range, irreducible and have mean velocity $v$:
\begin{equation*}
\sum_y y \, p(y,v)\;=\; v\; .
\end{equation*}
For the sake of concreteness, we take in this paper $p(y, v)=0$
unless $|y|=\sum_{1\le j\le d} |y_j|=1$ in which case
$$
p(\pm e_j,v)\;=\; \gamma \;\pm \; (1/2) e_j\cdot v
$$
for each vector $e_j$ and some $\gamma >0$ large enough for all
rates $p(e,v)$ to be non-negative.  The generator ${\mc L}^{ex}$ of
the random walk part of the dynamics acts on local functions $f$ on
the configuration space $X$ by
\begin{equation*}
(\mc L^{ex} f) (\eta) \;=\; \sum_{\substack{v \in \mc V \\
e \in \mc E}} \sum_{x \in \bb Z^d}
p(x, e, v; \eta)\, [f(\eta^{x, x+e, v}) - f(\eta)]\;,
\end{equation*}
where
\begin{equation*}
p(x, e, v;\eta)\;=\; \eta(x,v)\, [1-\eta(x+e,v)]\, p(e,v)
\end{equation*}
is the jump rate from $x$ to $x+e$ for particles with velocity $v$ and
\begin{equation*}
\eta^{x,y,v}(z,w) \;=\;
\left\{
\begin{array}{ll}
\eta(y,v) & \text{if $w=v$ and $z=x$,} \\
\eta(x,v) & \text{if $w=v$ and $z=y$,} \\
\eta(z,w) & \text{otherwise.}
\end{array}
\right.
\end{equation*}

The collision part of the dynamics is described as follows. Denote by
${\mc Q}$ a collision set which preserves momentum:
\begin{equation*}
\mc Q \;\subset\; \{(v,w,v',w') \in \mc V^4: v+w = v'+w'\}\;.
\end{equation*}
Assume that $\mc Q$ is symmetric in the sense that $(v,w,w',v')$,
$(v',w',v,w)$, and $(v',w',w,v)$ belong to $\mc Q$ as soon as
$(v,w,v',w')$ belongs to $\mc Q$.  Particles of velocities $v$ and $w$
at the same site collide at rate one and produce two particles of
velocities $v'$ and $w'$ at that site.  The generator $\mc L^c$ is
therefore
\begin{equation*}
\mc L^c f(\eta) \;=\; \sum_{y \in \bb Z^d}
\sum_{q \in \mc Q} p(y,q,\eta) \, [f(\eta^{y,q}) - f(\eta)]\;,
\end{equation*}
where the rate $p(y,q,\eta)$ is given by
\begin{equation*}
p(y,q,\eta) \;=\; \eta(y,v) \, \eta(y,w)\,
[1- \eta(y,v')]\, [1-\eta(y,w')]
\end{equation*}
and, for $q= (v_0, v_1, v_2, v_3)$, the configuration $\eta^{y,q}$
after the collision is defined as
\begin{equation*}
\eta^{y,q}(z,u) \;=\; \left\{
\begin{array}{ll}
\eta(y,v_{j+2})  & \text{if $z=y$ and $u=v_j$ for some $0\le j\le 3$,}\\
\eta(z,u) & \text{otherwise,}
\end{array}
\right.
\end{equation*}
where the sum in $v_{j+2}$ should be understood modulo $4$. 

The generator $\mc L$ of the lattice gas we examine in this article
is the superposition of the exclusion dynamics with the collisions
just introduced:
\begin{equation*}
\mc L  \;=\; \mc L^{ex} \;+\; \mc L^c \;.
\end{equation*}
Let $\{\eta (t) : t\ge 0\}$ be the Markov process with generator $\mc
L$.

\subsection{The invariant states}
We assume that the sets $\mc V$ and $\mc Q$ are chosen in such a way
that the unique conserved quantities are the local mass $I_0$ and
local momentum $I_a$, $a =1, \dots, d$:
\begin{equation*}
I_0(\eta_x) \;=\; \sum_{v\in \mc V} \eta(x,v)\;, \quad
I_{a}(\eta_x) \;=\; \sum_{v\in \mc V} (v\cdot e_{a})\,
\eta (x,v)\;.
\end{equation*}
Examples of sets $\mc V$ and collision dynamics with this property are
easy to produce.  Consider, for example, $\mc V = \mc E$ and take $\mc
Q$ as the set of all vectors $(v,w,v',w')$ such that $v+w=v'+w'=0$. An
elementary computation shows that the unique conserved quantities are
total mass and momentum.  \cite{EMYa} presents another example in
$d=3$.

For each chemical potential $\bs \lambda = (\lambda_0, \lambda_1,
\dots, \lambda_d)$, denote by $m_{\bs \lambda}$ the measure on
$\{0,1\}^{\mc V}$ given by
\begin{equation}
\label{refea}
m_{\bs \lambda} (\xi) \;=\; \frac 1{Z(\bs \lambda)}
\exp\Big\{ \sum_{a =0}^d \lambda_a I_{a}(\xi)
\Big\}\;, 
\end{equation}
where $Z(\bs \lambda)$ is a normalizing constant, $I_0(\xi) =
\sum_{v\in\mc V} \xi(v)$, $I_{a}(\xi) = \sum_{v\in\mc V} (v\cdot
e_{a}) \xi(v)$ for $a =1, \dots, d$. Notice that $m_{\bs
  \lambda}$ is a product measure on $\{0,1\}^{\mc V}$, i.e., that the
variables $\{\xi(v)\,: v\in \mc V\}$ are independent under $m_{\bs
  \lambda}$.

Denote by $\mu_{\bs \lambda}$ the product measure on $(\{0,1\}^{\mc
  V})^{\bb Z^d}$ with marginals given by
\begin{equation*}
\mu_{\bs \lambda} \{ \eta : \eta (x,\cdot) = \xi\}
\;=\; m_{\bs \lambda} (\xi)
\end{equation*}
for each $\xi$ in $\{0,1\}^{\mc V}$ and $x$ in $\bb Z^d$. $\mu_{\bs
  \lambda}$ is a product measure in the sense that the variables
$\{\eta (x,v)\, : x\in \bb Z^d, v\in \mc V\}$ are independent under
$\mu_{\bs \lambda}$. A simple computation shows that $\mu_{\bs
  \lambda}$ is an invariant state for the Markov process with
generator $\mc L$ for each $\bs \lambda$ in $\bb R_+ \times \bb R^d$,
that the generator $\mc L^c$ is symmetric with respect to $\mu_{\bs
  \lambda}$ and that $\mc L^{ex}$ has an adjoint $\mc L^{ex, *}$ in
which $p(e,v)$ is replaced by $p^*(e,v) = p(-e,v)$.

The expected value of the density under an invariant state can be
computed explicitly. Fix a vector $v$ in $\mc V$ and define $\theta_v
: \bb R^{d+1} \to \bb R_+$ as the expected value of $\eta(x,v)$ under
$\mu_{\bs \lambda}$:
\begin{equation*}
\theta_v (\bs \lambda) \;:=\; E_{\mu_{\bs \lambda}} [\eta(x,v)]
\;=\; \frac{\exp\big\{ \lambda_0 + \sum_{a =1}^d
\lambda_a (v\cdot e_a)\big\}}
{1+\exp\big\{ \lambda_0 + \sum_{a =1}^d
\lambda_a (v\cdot e_a)\big\} }\;\cdot
\end{equation*}
In this formula and below, for a probability measure $\mu$, $E_{\mu}$
stands for the expectation with respect to $\mu$. The expectation
under the invariant state $\mu_{\bs \lambda}$ of the mass and momentum
are given by
\begin{eqnarray*}
\rho(\bs \lambda) & := & E_{\mu_{\bs \lambda}}
[ I_0(\eta_x)] \;=\; \sum_{v\in {\mc V}}
\theta_v (\bs \lambda)  \\
u_a(\bs \lambda) & := & E_{\mu_{\bs \lambda}}
[ I_a (\eta_x) ] \;=\; \sum_{v\in {\mc V}} (v\cdot e_a )
\theta_v (\bs \lambda)\;.
\end{eqnarray*}

$(\rho(\bs \lambda) , u_1(\bs \lambda), \dots, u_d(\bs \lambda))$ is
the gradient of the strictly convex function $\log Z(\bs \lambda)$. In
particular, the map which associates the chemical potential $\bs
\lambda$ to the vector of density and momentum $(\rho, \bs u) = (\rho,
u_1, \dots, u_d)$ is one to one. Therefore, the chemical potential
$\bs \lambda = (\lambda_0, \dots, \lambda_d)$ can be expressed in
terms of $(\rho, \bs u)$: there exist a subset $\mf A$ of $\bb R_+
\times \bb R^{d}$ and functions $\Lambda_a : \mf A \to \bb R$, $a =0,
\dots, d$, such that
\begin{equation}
\label{eq:1}
\lambda_a \;=\; \Lambda_a (\rho (\bs \lambda),
\bs u (\bs \lambda))
\end{equation}
for each $\bs \lambda$ in $\bb R_+ \times \bb R^d$.  Let $\bs \Lambda
= (\Lambda_0, \dots, \Lambda_d)$.  This correspondence permits to
parameterize the invariant states by the density and the momentum: For
each $(\rho, \bs u)$ in $\mf A$, let 
$$
\nu_{\rho, \bs u} \;=\;
\mu_{\bs \Lambda(\rho, \bs u)}\;.  
$$

\subsection{Hydrodynamical limit under Euler scaling}
In this subsection, we deduce the hydrodynamic equation of the system
under the assumption of conservation of local equilibrium.  Fix smooth
functions $\rho : \bb R^d\to \bb R_+$, $\bs u : \bb R^d\to\bb R^d$.
For each $\varepsilon >0$, denote by $\nu^\varepsilon_{\rho(\cdot),
  \bs u(\cdot)}$ the product measure on $X$ with marginals given by
$$
\nu^\varepsilon_{\rho(\cdot), \bs u(\cdot)} \{ \eta (x, \cdot) = \xi\}
\;=\; \nu_{\rho(\varepsilon x), \bs u(\varepsilon x)} \{ \eta (0, \cdot) =
\xi\}
$$
for each $x$ in $\bb Z^d$ and each configuration $\xi$ in
$\{0,1\}^{\mc V}$. Assume that particles are initially distributed
according to $\nu^\varepsilon_{\rho(\cdot), \bs u(\cdot)}$.

For $j =1, \dots, d$, let $\nabla^-_j$ denote the lattice
gradient acting on functions $f:\bb Z^d\to\bb R$ by $\nabla^-_j
f(z)=f(z)-f(z-e_j)$. From Ito's formula, for $a=0,\dots ,d$,
we have the conservation law
\begin{equation*}
d I_a (\eta_x(t)) \;=\; \sum_{j =1}^d \nabla_j^-
w^a_{x,j} \, dt \;+\; dM^a_{x}(t)\;,
\end{equation*}
where $M^a_{x}(t)$ are martingales and $w^{a}_{x,j}$ are the
currents defined by
\begin{equation*}
\mc L \, I_{a}(\eta_x) \;=\;
\mc L^{ex} \, I_{a}(x , \eta) \;=\;
\sum_{j = 1}^d \nabla^-_j
\omega_{x,j}^{a}\;,
\end{equation*}
where
\begin{eqnarray*}
\omega^0_{x,j} &=& \gamma\, \nabla^-_{j} I_0(\eta_{x+e_j})
\;+\; \sum_{v\in \mc V} (e_{j}\cdot v) \, b_{x,j}(v) \;, \\
\omega^{a}_{x,j} &=& \gamma \, \nabla^-_{j} I_{a}
(\eta_{x+e_j}) \;+\;
\sum_{v\in \mc V} (e_{a}\cdot v)(e_{j}\cdot v)b_{x,j}(v)
\quad\text{and}\\
b_{x,j}(v) &=& \eta(x+e_{j},v) \eta(x,v)
\;-\;  (1/2) [\eta(x+e_{j},v) + \eta(x,v) ] \;.
\end{eqnarray*}
In this computation, the full generator can be replaced by the
exclusion operator because the collision operator preserves the
density and the momentum. Let $\omega^{a}_{j} = \omega^{a}_{0,j}$ for
$0\le a\le d$, $1\le j\le d$.

The expectation of the mass current in the $j$-th direction under the
local Gibbs state $\nu^\varepsilon_{\rho(\cdot), \bs u(\cdot)}$ is
denoted by $\pi_{0,j}$ and given by
\begin{eqnarray*}
\pi_{0,j} (\rho (\varepsilon x), \bs u (\varepsilon x)) 
&:=& E_{\nu^\varepsilon_{\rho(\cdot), \bs u(\cdot)}}
[ w^0_{x, j} ] \\
&=& \sum_{v\in \mc V} (v\cdot e_j )\,
\theta_v (\bs \Lambda (\rho, \bs u)) \,
\big\{\theta_v (\bs \Lambda (\rho, \bs u)) - 1 \big\}
\;+\; \gamma (\nabla_j^- \rho) (\varepsilon x)\; ,
\end{eqnarray*}
while the expectation of the momentum currents $\pi_{a,j}$
are given by
\begin{eqnarray*}
\!\!\!\!\!\!\!\!\!\!\!\! &&
\pi_{a ,j}(\rho (\varepsilon x), \bs u (\varepsilon x)) \; :=\;
E_{\nu^\varepsilon_{\rho(\cdot), \bs u(\cdot)}}[ w^a_{x, j} ] \\
\!\!\!\!\!\!\!\!\!\!\!\! && \qquad
=\; \sum_{v \in {\mc V}} (e_a \cdot v) \, (e_j \cdot v)
\, \theta_v (\bs \Lambda (\rho, \bs u)) \,
\big \{\theta_v (\bs \Lambda (\rho, \bs u)) - 1 \big \}
\;+ \; \gamma (\nabla_j^- u_a) (\varepsilon x)  \;.
\end{eqnarray*}
In both formulas, on the right hand side, $\rho$ and $\bs u$ are
evaluated at $\varepsilon x$.

Assuming conservation of local equilibrium, it is not difficult to
derive the hydrodynamic equations in the Euler scale for the lattice
gas considered in this article (cf. \cite{kl}). It is given by the
system of hyperbolic equations
\begin{equation*}
\left\{
\begin{array}{l}
\partial_t \rho \;+\; \sum_{j =1}^d \partial_{x_j} \pi_{0,j}
\;=\; 0\; , \\
\partial_t u_a \;+\; \sum_{j=1}^d \partial_{x_j}
\pi_{a, j } \;=\; 0\;.
\end{array}
\right.
\end{equation*}

Notice that the factors $\gamma \nabla_j^- $ do not survive in the
limit due to the presence of a second derivative. They will appear,
however, in the diffusive scale.

\subsection{Incompressible limit}
Inserting the local equilibrium assumption into the conservation laws,
on the time scale $\varepsilon^{-1} t$ one would obtain the equations
\begin{equation*}
\left\{
\begin{array}{l}
\partial_t \rho \;+\; \varepsilon^{-1} \sum_{j =1}^d
\partial_{x_j} \pi_{0,j} \;=\; \gamma \Delta \rho \; , \\
\partial_t u_a \;+\; \varepsilon^{-1} \sum_{j=1}^d
\partial_{x_j} \pi_{a, j } \;=\; \gamma \Delta  u_a \;.
\end{array}
\right.
\end{equation*}

Consider a system in which the density is near a constant and velocity
is of order $\varepsilon$. Expanding the mass density, the momentum
density and the mass and momentum currents, we obtain
\begin{align*}
\rho  & = \rho^{(0)} + \varepsilon \rho^{(1)} + \varepsilon^2
\rho^{(2)} +\cdots, \\
u  & =  \varepsilon u^{(1)} + \varepsilon^2 u^{(2)}+\cdots, \cr
\pi & = \pi^{(0)} + \varepsilon \pi^{(1)} + \varepsilon^2 \pi^{(2)}
+\cdots\;.
\end{align*}
Using these expansions and assuming that the zeroth order terms of
$\pi$ are constants, we obtain the incompressible Navier-Stokes
equations
\begin{equation*}
\left\{
\begin{array}{l}
\sum_{j=1}^d \partial_{x_j} \pi_{a, j}^{(1)} = 0  \\
\partial_t u^{(1)}_a
\;+\;  \sum_{j =1}^d \partial_{x_j} \pi^{(2)}_{a,j }
= \sum_{i , j}\sum_b   D^{a, b}_{i,j} \partial_{x_i}
\partial_{x_j} u^{(1)}_j
\end{array}
\right.
\end{equation*}
for $a=0,\dots ,d$ and $u_{0} =\rho$. Here the diffusion coefficient $
D^{a, b}_{i,j}= \gamma \delta_{a, b} \delta_{i, j}$ is a diagonal
matrix.  It turns out that this naive computation is correct if the
diffusion coefficient is instead given by a Green-Kubo formula.

\subsection{Green-Kubo formula.}
For simplicity, let $\lambda_j=0$ for $0\le j\le d$ so that $\theta_v
(\bs \lambda)= 1/2$ for every $v \in \mc V$. Denote this measure by
$\mu_0$ and let $\xi (x,v) = \eta(x,v) - \theta_v (\bs \lambda) =
\eta(x,v) - 1/2$.

Let $\tau_x$ be the shift by $x$ on the lattice, so that $(\tau_x
\eta) (y,v) = \eta(x+y,v)$. Denote by $\ll \cdot, \cdot\gg = \ll
\cdot, \cdot\gg_{\mu_0}$ the scalar product defined on $X$ by
$$
\ll f, g\gg \;=\; \sum_{x\in\bb Z^d} E_{\mu_0}[f ; \tau_x g]
$$
for two local functions $f$, $g$.  Here $E_{\mu_0}[g; h]$ stands
for the covariance between $g$ and $h$: $E_{\mu_0}[g; h] =
E_{\mu_0}[g\, h]-E_{\mu_0} [g ]E_{\mu_0}[h]$. Let ${\mc G}_0$ be the
space of local functions satisfying
$$
E_{\mu_0}[g]\;=\; 0\quad \text{and} \quad \ll g , I_a \gg
\;=\; 0
$$
for $0\le a\le d$.

Let $\chi = \{\chi_{a, b} , \, 0\le a,b\le d\}$ be the susceptibility
which in our context is given by
\begin{equation*}
\chi_{a, b} \;=\; \ll I_a , I_b \gg
\;=\; E_{\mu_0}  [ I_a(\eta_0)  ; I_b (\eta_0)]\;.
\end{equation*}
Denote by $\sigma^a_{j}$ the part of the current orthogonal
to the constants of motion:
\begin{equation*}
\sigma^a_{j} \;=\; \omega^a_{j}
\;-\; \sum_{b=0}^d c^{a,b}_{j} I_b (\eta_0)\;,
\end{equation*}
where the coefficients $c^{a,b}_{j}$ are chosen for $\sigma^a_{j}$ to
belong to ${\mc G}_0$. An elementary computation shows that
$$
c^{a,b}_j \;=\; \sum_{e =0}^d \ll \omega^a_{j} ,
I_e \gg (\chi^{-1})_{e, b} \;=\;
\frac{\partial}{\partial \alpha_b} E_{\nu_{\bs\alpha}}
\big[\omega^a_{j}\big]\; ,
$$
$\bs\alpha = (\alpha_0, \dots, \alpha_d)$.  Moreover, for $1\le
a\le d$,
\begin{eqnarray*}
\!\!\!\!\!\!\!\!\!\!\!\!\!\! &&
\sigma^0_{j} \;=\; \gamma \{ I_0(\eta_{e_j}) - I_0(\eta_0) \}
\;+\; \sum_{v\in \mc V} (e_{j} \cdot v) \xi(e_j,v) \xi(0,v)\;, \\
\!\!\!\!\!\!\!\!\!\!\!\!\!\! && \quad
\sigma^a_{j} \;=\; \gamma \{ I_a(\eta_{e_j}) - I_a(\eta_0) \}
\;+\; \sum_{v\in \mc V} (e_{j} \cdot v) (e_a \cdot v) \xi(e_j,v)
\xi(0,v)\;.
\end{eqnarray*}

For simplicity assume that the susceptibility is a constant times the
identity: $\chi_{a, b} = \kappa \, \delta_{a,b}$. A
straightforward computation shows that this is the case if the set $\mc
V$ is a cube centered at the origin. Under this assumption, following the
computation presented in section 2 of \cite{loy}, we obtain that for
$0\le a,b\le d$, $1\le i,j\le d$,
\begin{eqnarray*}
D_{i,j}^{a,b} (t) & := & \frac 1{2t\kappa}
\Big\{ \sum_{x\in\bb Z^d} E_{\mu_0} \big[ I_a (
\eta_x(t)) ; I_b(\eta_0(0)) \big] - t^2 V^{a,b}_{i,j} \Big\} \\
&=& \gamma \delta_{a,b} \delta_{i,j} + \frac 1{2t\kappa} \int_0^t
ds\int_0^s dr \ll \sigma^a_{i} , e^{r \mc L} \sigma^b_{j} \gg \\
&& \qquad\qquad  +\; \frac 1{2t\kappa} \int_0^t ds\int_0^s
dr \ll \sigma^a_{j} , e^{r \mc L}  \sigma^b_{i} \gg\;,
\end{eqnarray*}
where 
$$
2 V^{a,b}_{i,j} \;=\; \nabla E_{\nu_{\bs\alpha}} [
\omega^a_{0,i}] \cdot \nabla E_{\nu_{\bs\alpha}} [\omega^b_{0,j}]
\;+\; \nabla E_{\nu_{\bs\alpha}} [ \omega^b_{0,i}] \cdot \nabla
E_{\nu_{\bs\alpha}} [\omega^a_{0,j}] 
$$
and, for a local function $h$,
$\nabla E_{\nu_{\bs\alpha}} [h] = ((\partial/\partial \alpha_0)
E_{\nu_{\bs\alpha}} [h] , \dots, (\partial/\partial\alpha_d)
E_{\nu_{\bs\alpha}} [h])$. In dimension $d\ge 3$, the diffusion
coefficients $D_{i,j}^{a,b} (t)$ converge, as $t\uparrow\infty$, to
the diffusion coefficients $D_{i,j}^{a,b}$ given by the incompressible
Navier-Stokes equations in Subsection 2.3 (cf. \cite{EMYa}).

For $\theta$ in $\bb R^d$ and $r$ in $\bb R^{d+1}$, let
$$
D_{\theta, r} (t) \;=\; \sum_{a,b=0}^d \sum_{i,j=1}^d r_a \theta_i
D_{i,j}^{a,b} (t) \theta_j r_b\;.
$$
We can now state the main result. Let $\bb R_*^n = \bb R^n \setminus
\{0\}$.

\begin{theorem}
\label{s0}
Fix $\theta$ in $\bb R^d_*$ and $r$ in $\bb R^{d+1}_*$.  In dimension
$d=2$, there exists a positive constant $C_0 = C_0(\theta, r)$ so that
for all sufficiently small $\lambda>0$,
\begin{equation*}
\int_0^\infty dt \, e^{-\lambda t} \; t D_{\theta,r} (t)
\; \ge\;  C_0 \lambda^{-2} \log \log \lambda^{-1}  \; .
\end{equation*}
\end{theorem}

Recall that $\int_0^\infty e^{-\lambda t} f(t) dt \sim
\lambda^{-\alpha}$ as $\lambda\to 0$ means, in some weak sense, that
$f(t)\sim t^{\alpha -1}$ as $t\to \infty$. Theorem \ref{s0} is
therefore stating that the diffusion coefficient $D_{\theta, r} (t)$
is diverging as $\log \log t$ in a weak sense.

Fix a vector $\theta$ in $\bb R^d$, $r$ in $\bb R^{d+1}$ and let
$\sigma = \sigma_{\theta,r}$ be given by
\begin{equation}
\label{eq:11}
\sigma \;=\; \sum_{a=0}^d \sum_{j=1}^d r_a \theta_i \sigma^a_{i}\;.
\end{equation}
An elementary computation shows that for every $\lambda>0$
$$
\int_0^\infty dt \, e^{-\lambda t} \; t D_{\theta, r} (t)
\;=\; \gamma \lambda^{-2} \Vert \theta \Vert^2 \Vert r \Vert^2
\;+\; \lambda^{-2} \ll \sigma , (\lambda - \mc L)^{-1}
\sigma \gg\;.
$$
Therefore, Theorem \ref{s0} follows from next estimate on the
resolvent.

\begin{lemma}
\label{s-1}
Fix $\theta$ in $\bb R^d_*$ and $r$ in $\bb R^{d+1}_*$.  There exists a
positive constant $C_0 = C_0(\theta, r)$ such that for sufficiently
small $\lambda>0$,
\begin{equation*}
\ll \sigma ,(\lambda-{\mathcal L})^{-1} \sigma \gg
\; \ge\;  C_0  \, \log \log \lambda^{-1} \; .
\end{equation*}
\end{lemma}

Notice that the piece of the current $\gamma \{I_a(\eta_{e_j}) -
I_a(\eta_0)\}$ vanishes for the inner product $\ll \cdot, \cdot\gg$.
We may therefore ignore it in the computations below.

\section{Dual representation}
\label{sec2}

Recall that we are fixing the chemical potentials to be zero:
$\lambda_j=0$ for $0\le j\le d$ and that we denote by $\mu_0$ the
product invariant measure associated to this chemical potential.
Unless otherwise stated, $<f,g>$ stands for the inner product of $f$
and $g$ in $L^2(\mu_0)$.

Denote by $\mc E$ the finite subsets of $\bb Z^d \times \mc V$ and by
$\mc E_n$ the subsets of $\mc E$ with cardinality $n$, for $n\ge 0$.
For a set $A$ in $\mc E$, let $\Psi_A$ be the local function defined
by 
$$
\Psi_A (\eta) \;=\; \prod_{(x,v)\in A} [\eta (x,v) - 1/2] \;=\;
\prod_{(x,v)\in A} \xi (x,v) \;.  
$$
By convention $\Psi_\phi =1$,
where $\phi$ stands for the empty set.  It is not difficult to check
that $\{\Psi_A: A\in\mc E\}$ forms an orthogonal basis of
$L^2(\mu_0)$. In particular, any local function $f$ can be written as
$\sum_{A\in\mc E} \mf f(A) \Psi_A$ for some finite supported function
$\mf f : \mc E\to \bb R$. This latter function is called the Fourier
coefficients of the local function $f$ and frequently denoted by $\bb
T f$. A local function $f$ is said to have degree $n\ge 0$ if $(\bb T
f)(A) =0$ for any $A$ in $\mc E_n^c$.

Notice that for any local functions $f = \sum_{A\in\mc E} \mf f(A)
\Psi_A$, $g = \sum_{A\in\mc E} \mf g(A) \Psi_A$,
$$
<f,g>_{\mu_0} \;=\; \sum_{n\ge 0} (1/4)^n
\sum_{A\in\mc E_n} \mf f(A) \, \mf g(A)\; .
$$
The factor $(1/4)$ appears because we did not consider an
orthonormal basis since $E_{\mu_0}[ \xi(x,v)^2] = 1/4$.

We say that two finite subsets $A$, $B$ of $\bb Z^d \times \mc V$ are
equivalent if one is the translation of the other. This equivalence
relation is denoted by $\sim$ so that $A\sim B$ if $A=B+x$ for some
$x$ in $\bb Z^d$. Let $\tilde {\mc E}_n$ be the quotient of $\mc E_n$
with respect to this equivalence relation: $\tilde {\mc E}_n = \mc E_n
/_\sim$, $\tilde {\mc E} = {\mc E} /_\sim$.  An elementary computation
(cf. \cite{lov6}) gives that
$$
\ll f, g \gg_{\mu_0} \;=\; \sum_{n\ge 1} (1/4)^{n}
\sum_{A\in\tilde {\mc E}_n} \bar {\mf f} (A) \, \bar {\mf g} (A)\;,
$$
where
\begin{equation}
\label{refeb}
\bar {\mf f} (A) \;=\; \sum_{z\in\bb Z^d} \mf f (A+z)\;.
\end{equation}

Therefore, if we denote by $\ll \cdot, \cdot \gg$ the inner product in
$L^2(\mc E)$ defined by
\begin{equation*}
\ll \mf f, \mf g \gg \;=\; \sum_{n\ge 1} (1/4)^{n}
\sum_{A\in\tilde {\mc E}_n} \bar {\mf f} (A) \, \bar {\mf g} (A)\;,
\end{equation*}
where $\bar {\mf f}$, $\bar {\mf g}$ are defined by \eqref{refeb}, we
have that
\begin{equation*}
\ll f, g \gg_{\mu_0} \;=\; \ll  {\bb T f}, { \bb T g} \gg\;.
\end{equation*}

The goal of this section is to examine the action of the generators
$\mc L^{ex}$, $\mc L^c$ on the Fourier coefficients. More precisely,
to find operators $\bb L^{ex}$, $\bb L^c$ such that $\bb L^{ex} \bb T
= \bb T \mc L^{ex}$, $\bb L^{c} \bb T = \bb T \mc L^{c}$.

\subsection{The exclusion operator}
We start with the exclusion part of the generator which can be
decomposed into its symmetric part $\mc S$ and its antisymmetric part
$\mc A$ as given by
\begin{eqnarray*}
(\mc Sf)(\eta) &=& \gamma \sum_{v\in\mc V} \sum_{x\in\bb Z^d}
\sum_{j=1}^d  \{ f(\eta^{x,x+e_j,v})- f(\eta)\}\;, \\
(\mc Af)(\eta) &=& \frac 12 \sum_{v\in\mc V} \sum_{j=1}^d (e_j \cdot v)
\sum_{x\in\bb Z^d} \{\eta (x,v) - \eta (x+e_j,v)\} \,
\{ f(\eta^{x,x+e_j,v})- f(\eta)\}\;.
\end{eqnarray*}

To examine the action of the symmetric part of the exclusion generator
on the Fourier coefficients, we first compute $\mc S \Psi_A$.  An
elementary computation shows that for each set $A$ in $\mc E$,
$$
\mc S \Psi_A \;=\; \gamma \sum_{v\in\mc V} \sum_{j=1}^d
\sum_{x\in\bb Z^d} \{ \Psi_{A_{x,x+e_j,v}} - \Psi_A\}
$$
provided $A_{x,x+e_j,v}$ stands for
\begin{equation}
\label{eq:3}
\left\{
\begin{array}{ll}
(A \setminus \{(x,v)\})\cup \{(x+e_j,v)\} & \text{ if $(x,v) \in
  A$ and $(x+e_j,v) \not \in A$,} \\
(A \setminus \{(x+e_j,v)\} )\cup \{(x,v)\} & \text{ if $(x+e_j,v) \in
  A$ and $(x,v) \not \in A$,} \\
A & \text{otherwise.}
\end{array}
\right.
\end{equation}
In particular, for any local function $f = \sum_A \mf f(A) \Psi_A$, a
change of variable $B = A_{x,x+e_j,v}$ gives that
$$
\mc S f \;=\; \gamma \sum_{A\in\mc E} \Psi_A \sum_{v\in\mc V}
\sum_{x\in\bb Z^d} \sum_{j=1}^d \{\mf f (A_{x,x+e_j,v}) - \mf f(A)\}\;.
$$
Therefore, if we define the operator $\bb S$ as
$$
(\bb S \mf f)(A) \;=\; \gamma \sum_{v\in\mc V} \sum_{x\in\bb Z^d}
\sum_{j=1}^d \{\mf f (A_{x,x+e_j,v}) - \mf f(A)\}\;,
$$
we have that $\bb T \mc S = \bb S \bb T$.

We turn now to the antisymmetric part.  To compute $\mc A \Psi_A
(\eta)$, observe that $[\eta(y,w) -1/2] \Psi_A$ is equal to
$\Psi_{A\cup\{(y,w)\}}$ if $(y,w)$ does not belong to $A$ and is equal
to $(1/4) \Psi_{A\setminus \{(y,w)\}}$ if $(y,w)$ belongs to $A$
because $(\eta(y,v) -1/2)^2 = 1/4$. In particular, a straightforward
computation shows that
\begin{eqnarray*}
\mc A \Psi_A &=&  \sum_{j=1}^d \sum_{\substack{(x,v)\in A \\
(x+e_j, v)\not \in A}} (e_j\cdot v) \{ \Psi_{A \cup \{(x+e_j,v)\}} -
(1/4) \Psi_{A \setminus \{(x,v)\}}\} \\
&+&  \sum_{j=1}^d \sum_{\substack{(x+e_j,v)\in A \\
(x, v)\not \in A}} (e_j\cdot v) \{ (1/4)  \Psi_{A \setminus \{(x+e_j,v)\}} -
\Psi_{A \cup \{(x,v)\}}\}\;.
\end{eqnarray*}
In the first term on the right hand side, the second sum is carried
over all pairs $(x,v)$ in $\bb Z^d\times \mc V$, such that $(x,v)$
belongs to $A$ and $(x+e_j,v)$ does not.  Therefore, if $f =
\sum_{A\in\mc E} \mf f(A) \Psi_A$ is a local function, after elementary
changes of variables, we obtain that
$$
\mc A f \;=\; \sum_{A\in\mc E} (\bb A \mf f)(A) \Psi_A
$$
provided $\bb A = \bb J_+ + \bb J_-$, where
\begin{eqnarray*}
\!\!\!\!\!\!\!\!\!\!\!\!\!\! &&
(\bb J_+ \mf f)(A) \;=\; \sum_{j=1}^d \sum_{\substack{(x,v)\in A \\
(x+e_j, v) \in A}} (e_j\cdot v) \big \{ \mf f (A \setminus \{(x+e_j,v)\}) -
\mf f(A \setminus \{(x,v)\}) \big \}\;, \\
\!\!\!\!\!\!\!\!\!\!\!\!\!\! && \quad
(\bb J_- \mf f)(A) \;=\; (1/4) \sum_{j=1}^d \sum_{\substack{
(x,v)\not \in A \\ (x+e_j, v) \not \in A}} (e_j\cdot v)
\big \{ \mf f (A \cup \{(x+e_j,v)\}) -
\mf f(A \cup \{(x,v)\}) \big \}
\end{eqnarray*}
and $\bb T \mc A = \bb A\bb T$.

Let $\bb L^{ex} = \bb S + \bb A$. Up to this point we proved that $\bb
T \mc L^{ex} = \bb L^{ex} \bb T$.  This operator $\bb L^{ex}$, which is
not a generator, can thus be decomposed in three pieces, $\bb S$, $\bb
J_+$, $\bb J_-$. $\bb S$ is the symmetric part of $\bb L^{ex}$ and
does not change the degree of a function. In contrast, $\bb J_+$
increases the degree by one, while $\bb J_-$ decreases it by one and
$\bb J_-$ is the adjoint of $- \bb J_+$: $(\bb J_-)^* = - \bb J_+$. In
particular, $\bb J_- + \bb J_+$ is anti-symmetric:
\begin{eqnarray*}
&& \ll \bb S \mf f, \mf g \gg \; =\; \ll \mf f, \bb S \mf g \gg
\quad\text{and}\quad \ll \bb J_+ \mf f, \mf g \gg \; =\; -
\ll \mf f, \bb J_- \mf g \gg \\
&&\quad\text{so that}\quad \ll (\bb J_+ + \bb J_-) \mf f, \mf g \gg \;
=\; - \ll \mf f, (\bb J_+ + \bb J_-) \mf g \gg\; .
\end{eqnarray*}
Moreover, a simple computation shows that in $L^2(\mc E)$
\begin{eqnarray}
\label{eq:8}
&& \bb S \mf f \; =\; 0 \quad\text{for all functions $\mf f$ of degree
  one}\; , \\
&& \quad \bb J_- \mf f \; =\; 0 \quad\text{for all functions $\mf f$
of degree two}\; .
\nonumber
\end{eqnarray}

\subsection{The collision operator}
The remainder of this section is devoted to the collision operator.
We start defining a generator $\mc L_1^c$ and showing in Lemma
\ref{s1} below that it is of the same order as $\mc L^c$ for our
purposes. We conclude the section investigating the action of $\mc
L_1^c$ on the Fourier coefficients.

Fix a site $x$ in $\bb Z^d$ and a point $q= (v,w,v',w')$ in the
collision set $\mc Q$. Let $q' = (v,w,w',v')$. Denote by $\mc L_{x,q}$
the generator defined by
$$
(\mc L_{x,q} f)(\eta) \;=\; p'(x,q,\eta) \Big\{ [f(\eta^{x,q}) -
f(\eta)] + [f(\eta^{x,q'}) - f(\eta)] \Big\}\;,
$$
where
\begin{eqnarray*}
p'(x,q,\eta) &=& \eta(x,v)\, \eta(x,w) \, [1-\eta(x,v')]
\, [1-\eta(x,w')] \\
&+& \eta(x,v')\, \eta(x,w') \, [1-\eta(x,v)] \, [1-\eta(x,w)]\;.
\end{eqnarray*}
Since the collision set $\mc Q$ is symmetric, we may rewrite the
collision generator $\mc L^c$ as
$$
\mc L^c \;=\; (1/4) \sum_{x\in\bb Z^d} \sum_{q\in Q} \mc L_{x,q}\; .
$$

Fix a site $x$ in $\bb Z^d$.  We start the analysis of the collision
operator by examining the generator $\mc L_{x,q}$. Since site $x$ is
fixed, we omit the index $x$ below so that $\mc L_q = \mc L_{x,q}$. We
also denote by $\zeta$ configurations of $\{0,1\}^{\mc V}$ and by
$\xi(v)$ the function $\zeta (v) - 1/2$.  Since $m_{\bs \lambda}$
given by \eqref{refea} is a product measure on $\{0,1\}^{\mc V}$, $\{
\prod_{v\in B} \xi (v)\, : B \subset V\}$ forms an orthogonal basis of
$L^2(m_0)$ if $\prod_{v\in \phi} \xi (v) = 1$, where $\phi$ stands for
the empty set.

Fix $q= (v,w,v',w')$ in $\mc Q$ and let $H_q =\{v,w,v',w'\}$. Since
$\mc L_q$ only changes the variables $\{\zeta(u) \, :u\in H_q\}$, all
other variables can be considered as constants so that for any subset
$B$ of $\mc V$,
$$
\mc L_q \prod_{u\in B} \zeta (u) \;=\; \prod_{u\in B\setminus H_q}
\zeta (u) \, \mc L_q \prod_{u\in B\cap H_q} \zeta (u)\;.
$$
A similar identity holds if we replace $\zeta$ by $\xi$. 

Let
\begin{eqnarray*}
\phi_1 (\xi) &=& \xi (v') \;+\; \xi(w') \;-\; \xi(v) \;-\; \xi(w)\;,
\\
\tilde \phi_3 (\xi) &=\;& \xi (v)\xi (w) \xi (v') \;+\;
\xi (v)\xi (w) \xi (w') \;-\; \xi (v')\xi (w') \xi (v)
\;-\; \xi (v')\xi (w') \xi(w)\; .
\end{eqnarray*}
The index $1$ and $3$ stand for the degree of the functions involved.
Straightforward computations give the following identities for degree
one functions:
\begin{eqnarray*}
\!\!\!\!\!\!\!\!\!\!\!\!\!\! &&
\mc L_q \xi (v) \;=\; \mc L_q \xi (w) \;=\; (1/2) \phi_1 (\xi) \;+\;
2 \tilde \phi_3 (\xi) \;, \\
\!\!\!\!\!\!\!\!\!\!\!\!\!\! && \quad
\mc L_q \xi (v') \;=\; \mc L_q \xi (w') \;=\; - (1/2) \phi_1 (\xi) \;-\;
2 \tilde \phi_3 (\xi)\;.
\end{eqnarray*}
Degree two functions vanish under the action of the generator:
$$
\mc L_q \xi(u_1) \xi(u_2) \;=\; 0
$$
for $u_1$, $u_2$ in $H_q$. To derive these identities we used
that $\mc L_q \{ \xi(v) + \xi(v')\} =0$ and similar equalities.
Degree three functions are such that
\begin{eqnarray*}
\!\!\!\!\!\!\!\!\!\!\!\!\!\! &&
\mc L_q \xi (v) \xi (w) \xi (v') \;=\; \mc L_q \xi (v) \xi (w) \xi
(w') \;=\; - (1/8) \phi_1 (\xi) \;-\; (1/2) \tilde \phi_3 (\xi) \;, \\
\!\!\!\!\!\!\!\!\!\!\!\!\!\! && \quad
\mc L_q \xi (v') \xi (w') \xi (v) \;=\; \mc L_q \xi (v') \xi (w') \xi
(w) \;=\; (1/8) \phi_1 (\xi) \;+\; (1/2) \tilde \phi_3 (\xi) \;.
\end{eqnarray*}
Finally, degree four functions vanish under the action of the generator:
$$
\mc L_q \xi(v) \xi(w) \xi(v') \xi(w') \;=\; 0\;.
$$
Here again, to deduce this equality we used that $\mc L_q \{ \xi
(v) \xi (w) \xi (v') + \xi (v') \xi (w') \xi (v) \}$ vanishes as well
as similar identities.

It follows from the previous formulas that the unique non zero
eigenvalue $-4$ is associated to the eigenfunction $\psi= \phi_1 +
4\tilde \phi_3$. In particular, the generator $\mc L_q$ can be written
as
$$
\mc L_q f \;=\; - 4 \frac{<f,\psi>}{<\psi,\psi>} \psi \;,
$$
where $<\cdot, \cdot>$ stands for the inner product in
$L^2(m_0)$. Denote by $\mc L_{q,1}$, $\mc L_{q,3}$ the operators
defined by
$$
\mc L_{q,1} f \;=\; - 4 \frac{<f,\phi_1>}{<\phi_1,\phi_1>} \phi_1
\;,\quad \mc L_{q,3} f \;=\; - 4 \frac{<f,\phi_3>}{<\phi_3,\phi_3>}
\phi_3
$$
for $\phi_3 = 4 \tilde \phi_3$.  Since $\phi_1$, $\phi_3$ are
orthogonal, an elementary computation shows that
\begin{equation}
\label{eq:4}
-\mc L_q \;\le\; -2 \mc L_{q,1} \;-\; 2 \mc L_{q,3}
\end{equation}
in the matrix sense.

If we reintroduce the index $x$, we obtain the operators
$$
\mc L^c_j \;=\; (1/4) \sum_{x\in\bb Z^d} \sum_{q\in Q} \mc L_{x,q,j}
$$
for $j=1$, $3$. Notice that both operators keep the degree of local
functions. Indeed, for a local function $f$, an elementary computation
shows that
$$
\mc L^c_j f \;=\; - \sum_{x\in\bb Z^d} \sum_{q\in Q} <f,
\phi_{x,q,j}>_x \phi_{x,q,j}
$$
because $<\phi_{x,q,j} ,\phi_{x,q,j}>_x =1$ for all $x$, $q$ and
$j$.  In this formula $<\cdot, \cdot>_x$ stands for the inner product
with respect to $m_0(\eta(x, \cdot))$, which means that only the
variables at site $x$ are integrated and $\phi_{x,q,j}$ are the
functions defined by
\begin{eqnarray*}
\phi_{x,q,1} (\eta) &=& \xi (x,v') \;+\; \xi(x,w')
\;-\; \xi(x,v) \;-\; \xi(x,w)\;, \\
\phi_{x,q,3} (\eta) &=& 4 \xi (x,v)\xi (x,w) \xi (x,v') \;+\;
4 \xi (x,v)\xi (x,w) \xi (x,w') \\
&-& 4 \xi (x,v')\xi (x,w') \xi (x,v)
\;-\; 4 \xi (x,v')\xi (x,w') \xi(x,w)\; .
\end{eqnarray*}
If $f = \Psi_B$ for some finite set $B$, an elementary computation
shows that
$$
<\Psi_B , \phi_{x,q,j}>_x \phi_{x,q,j}
$$
has the same degree as $\Psi_B$ for all $x$, $q$ and $j$, which
proves the claim.

It follows from (\ref{eq:4}) that
$$
- \mc L^c \;\le\; - 2 \mc L^c_{1} \;-\; 2 \mc L^c_{3}\;.
$$

In order to have a tractable algebraic expression for the collision
operator, we plan to substitute $\mc L^c$ by $\mc L^c_1$. In order to
estimate the third degree terms we use the following lemma.

\begin{lemma}
\label{s1}
There exists a finite constant $C_0$ such that
\begin{eqnarray*}
\!\!\!\!\!\!\!\!\!\!\!\!\! &&
{C_0}^{-1} \, \ll f , \{\lambda - \mc L^{ex}  -\mc L^c_1 \}^{-1}
f \gg \\
\!\!\!\!\!\!\!\!\!\!\!\!\! && \qquad\qquad\qquad
\leq\;\;  \ll f, \{\lambda - \mc L\}^{-1} f\gg \;\; \leq\;\;
C_0 \, \ll f , \{\lambda - \mc L^{ex}  -\mc L^c_1 \}^{-1} f \gg
\end{eqnarray*}
for every $\lambda>0$ and every local function $f$.
\end{lemma}

The proof of this result is based on the lemma below whose proof is
similar to the one of Lemma 4.2 in \cite{LQSY}

\begin{lemma}
\label{s2}
Consider a function $\omega: \mc E_3\to\bb R$.  Assume that there
exists $\ell \ge 1$ such that $\omega ((x_1,v_1), (x_2,v_2),
(x_3,v_3))=0$ if $|x_1-x_2| > \ell$ or $|x_1 -x_3|>\ell$. Then, there
exists a finite constant $C_0$ depending only on $\omega$ such that
\begin{equation*}
\sum_{(\bs x, \bs v)\in \mc E_3} \mf f^2(\bs x, \bs v)\,
\omega(\bs x, \bs v) \;\leq\;  C_0
\sum_{(\bs x, \bs v)\in \mc E_3} \mf f(\bs x, \bs v) (-\bb S \mf f)
(\bs x, \bs v)
\end{equation*}
for every finite supported function $\mf f: \mc E_3\to\bb R$.
\end{lemma}

\noindent{\bf Proof of Lemma \ref{s1}.}
We claim that $-\mc L^c_3 \le - C_0 \mc S$ for some finite constant
$C_0$. Since $\mc L^c_3$ keeps the degree, to prove the claim we only
need to show that
$$
<-\mc L_3^c f, f> \;\le\; C_0 <-\mc S f, f>
$$
for local functions of a fixed degree.

Fix $n\ge 1$ and a local function $f$ of degree $n$. By definition of
$\mc L^c_3$, taking conditional expectations we obtain that
$$
<-\mc L_3^c f, f> \;=\; \sum_{q\in \mc Q} \sum_{x\in\bb Z^d}
E\Big[ < f, \phi_{x,q,3}>_x^2 \Big]\;.
$$
Since the velocity set is finite, by definition of $\phi_{x,q,3}$,
to prove the claim it is enough to show that
$$
\sum_{x\in\bb Z^d} E\Big[ < f, \Psi_{B_x}>_x^2 \Big] \;\le\;
C_0 < (-\mc S)f,f>\;,
$$
where $B_x=\{(x,v_1), (x,v_2), (x,v_3)\}$ and $v_1$, $v_2$, $v_3$
are three distinct velocities in $\mc V$.  Assume that $f = \sum_B \mf
f(B) \Psi_B$. An elementary computation shows that the expectation
appearing on the left hand side of the previous inequality is bounded
above by
$$
(1/4)^n \sum_{B\supset B_x} \mf f(B)^2\;,
$$
where the sum is performed over all sets $B$ which contain $B_x$.
In particular, if for a finite set $A = \{(x_1, v_1),
(x_2, v_2), (x_3, v_3) \}$, we set
$$
\rho(A)^2 \;=\; (1/4)^n \sum_{B\supset A} \mf f(B)^2\;,
$$
where the summation is performed over all sets $B$ which contain $A$,
we just proved that
$$
<-\mc L_3^c f, f> \;\le\; \sum_{v_1, v_2, v_3} \sum_{x\in\bb Z^d}
\rho( \{(x,v_1), (x,v_2), (x,v_3)\})^2\;.
$$
By Lemma \ref{s2}, this expression is less than or equal to
$$
C_0 \sum_{(\bs x, \bs v)\in \mc E_3} \rho (\bs x, \bs v) (-\bb S \rho)
(\bs x, \bs v) \;\le\; C_0 <-\mc S f, f>\;.
$$
Last inequality follows from Schwarz inequality and concludes the
proof of the claim. Here and below, $C_0$, $a_0$ are constants whose
value may change from line to line.

Since $-\mc L^c_3 \le - C_0 \mc S$, we have that $-\mc L^c \le -2(\mc
L^c_1 + \mc L^c_3) \le -C_0 (\mc S + \mc L^c_1)$ and
$$
\lambda \;-\; \mc S \;-\; \mc L^c \;\le\;
a_0 \{ \lambda \;-\; \mc S \;-\; \mc L_1^c \}
$$
for some finite constant $a_0 >1$ and all $\lambda>0$.

On the other hand, since $\mc L^c_3 \ge  C_0 \mc S$ and since $(a-b)^2
\geq (1-\varepsilon) a^2 - (\varepsilon^{-1} -1) b^2$ for every $0<
\varepsilon <1$, a straightforward computation shows that
\begin{equation*}
-\mc L^c \;\geq\; - \frac{ (1-\varepsilon)}{2} \mc L^c_1
\; +\; \frac{ (\varepsilon^{-1} -1)}2 \mc L^c_3
\; \geq\; - \frac{ (1-\varepsilon)}{2} \mc L^c_1
\; +\; C_0 (\varepsilon^{-1} -1) \mc S
\end{equation*}
for some finite constant $C_0$. Here the factor $1/2$ appeared because
$<\psi, \psi>= 2$ $<\phi_1, \phi_1> = 2 <\phi_3, \phi_3>$. If we choose
$\varepsilon$ small enough for $C_0 (\varepsilon^{-1} -1) <1$, it
follows from this inequality that
\begin{eqnarray*}
\lambda \;-\; \mc S \;-\; \mc L^c &\geq& \lambda \;-\;
\frac{ (1-\varepsilon)}{2} \mc L^c_1 \;-\;
\big\{ 1- C_0 (\varepsilon^{-1} -1) \big\} \mc S \\
& \geq & a_0^{-1} \big\{ \lambda \;-\; \mc S \;-\;
\mc L^c_1  \big\}
\end{eqnarray*}
for some finite constant $a_0 >1$ and all $\lambda>0$.

Up to this point we proved the existence of a finite constant
$a_0>1$ such that
\begin{equation}
\label{eq:2}
a_0^{-1} \big\{ \lambda \;-\;  \mc L^c_1 \;-\;
\mc S \big\} \;\le\; \lambda \;-\; \mc S \;-\; \mc L^c \;\le\;
a_0 \{ \lambda \;-\; \mc S \;-\; \mc L_1^c \}
\end{equation}
for all $\lambda>0$.

It remains to add the asymmetric part of the exclusion generator.
Denote by $R^s$ (resp. $R^a$, $R^*$) the symmetric part (resp.
asymmetric part, adjoint) of an operator $R$.  It is well known that
$$
\big\{(R^{-1})^s\big\}^{-1} \; =\;  R^{*}  (R^s)^{-1} R\;=\;
R^s \;+\; (R^a)^* (R^s)^{-1} R^a\; .
$$
In particular, for every $\lambda>0$,
$$
\Big(\big\{(\lambda - \mc L)^{-1}\big\}^s\Big)^{-1}
\;=\; (\lambda - \mc S - \mc L^c) \; +\;
\mc A^*(\lambda - \mc S - \mc L^c)^{-1} \mc A\;.
$$
In view of (\ref{eq:2}), there exists a finite constant
$a_0>1$ such that
\begin{eqnarray*}
\!\!\!\!\!\!\!\!\!\!\!\!\! &&
a_0^{-1} \Big\{ (\lambda - \mc S - \mc L^c_1) \; +\;
\mc A^*(\lambda - \mc S - \mc L^c_1)^{-1} \mc A \Big\} \\
\!\!\!\!\!\!\!\!\!\!\!\!\! && \quad
\;\le\; (\lambda - \mc S - \mc L^c) \; +\;
\mc A^*(\lambda - \mc S - \mc L^c)^{-1} \mc A \\
\!\!\!\!\!\!\!\!\!\!\!\!\! && \qquad
\;\le\; a_0 \Big\{ (\lambda - \mc S - \mc L^c_1) \; +\;
\mc A^*(\lambda - \mc S - \mc L^c_1)^{-1} \mc A \Big\}
\end{eqnarray*}
so that
$$
a_0^{-1} \big\{(\lambda - \mc L^{ex} - \mc L^c_1)^{-1}\big\}^s \;\le\;
\big\{(\lambda - \mc L)^{-1}\big\}^s \;\le\;
a_0 \big\{(\lambda - \mc L^{ex} - \mc L^c_1 )^{-1}\big\}^s \;,
$$
which proves the lemma in the case of the inner product of $L^2(\mu_0)$.
The extension to the inner product $\ll \cdot, \cdot \gg$ is standard
(cf. \cite{LY}). \qed
\smallskip

We conclude the section examining the action of $\mc L^c_1$ on the
Fourier coefficients. For any local function $f = \sum_{A\in \mc E}
\mf f(A) \Psi_A$, a simple computation shows that
$$
\mc L^c_1 f \;=\; \sum_{A\in \mc E} (\bb L^c_1 \mf f)(A)
\Psi_A \;,
$$
where
\begin{eqnarray*}
\!\!\!\!\!\!\!\!\!\!\!\!\!\! &&
(\bb L^c_1 \mf f)(A) \;=\; (1/4) \sum_{q\in \mc Q} \sum_{x}
i_q(A_x) \big \{ \mf f(A_x^c \cup \{(x,v')\}) + \mf f(A_x^c \cup
\{(x,w')\}) \\
\!\!\!\!\!\!\!\!\!\!\!\!\!\! &&
\qquad\qquad\qquad\qquad\qquad\qquad\qquad\qquad\qquad
- \mf f(A_x^c \cup \{(x,v)\})- \mf f(A_x^c \cup \{(x,w)\})
\big \}\;.
\end{eqnarray*}
In this formula, $q=(v,w,v',w')$, $A_x$ stands for the set of
velocities $u$ such that $(x,u)$ belongs to $A$: $A_x =\{u\in\mc V:
(x,u)\in A\}$ , $A_x^c$ for the set of points $(y,u)$ in $A$ with
$y\neq x$: $A_x^c =\{(y,u) \in A: y\neq x\}$ and if $i_q(A_x)=1$
if $A_x$ is an incoming velocity, $-1$ is $A_x$ is an outgoing
velocity and $0$ otherwise:
$$
i_q(A_x) \;=\; \mb 1\{ A_x = \{v\}\}
\;+\; \mb 1\{ A_x = \{w\}\}
\;-\; \mb 1\{ A_x = \{v'\}\}
\;-\; \mb 1\{ A_x = \{w'\}\} \;.
$$
With this notation we have that $\bb T \mc L^c_1 = \bb L^c_1 \bb T$.
Notice that $-\bb L_1^c$ is a non-negative symmetric operator in
$L^2(\mc E_n)$ and that
$$
< -\bb L_1^c \mf f, \mf f> \;=\; \sum_x \sum_q E_{\mu_0}
\Big[ < f, \phi_{x,q,1}>^2_x \Big]\;.
$$

\section{Cutoff of large degrees}
\label{sec3}

For $n\ge 1$, let $\mc G_n = \cup_{1\le k\le n} \mc E_n$. Denote by
$\pi_n$ the orthogonal projection on $L^2(\mc E_n)$, by $P_n$ the
orthogonal projection on $L^2(\mc G_n)$ and by $\mc L_n$ the operator
$\mc L^{ex} + \mc L_1^c$ truncated at level $n$: $\mc L_n = P_n (\mc
L^{ex} + \mc L_1^c) P_n$. In particular, $\mf f = \sum_{n\ge 1} \pi_n
\mf f$ and $P_n = \sum_{1\le j\le n} \pi_j$.

To investigate the asymptotic behavior of $\ll \sigma, (\lambda -\mc
L^{ex} - \mc L_1^c)^{-1} \sigma \gg$, for $\lambda>0$ consider the
resolvent equation $(\lambda-\mc L^{ex} - \mc L_1^c) u_\lambda =
\sigma$.  In the Fourier space, the equation becomes the hierarchy
equations
\begin{eqnarray*}
\left\{
\begin{array}{l}
\vphantom{\Big\{}
\bb J_+^* \pi_3 \mf u_\lambda + (\lambda-\bb S - \bb L_1^c)
\pi_2 \mf u_\lambda = \sigma \; , \\
\vphantom{\Big\{}
\bb J_+^* \pi_{k+1} \mf u_\lambda + (\lambda -\bb S - \bb L_1^c)
\pi_k \mf u_\lambda - \bb J_+ \pi_{k-1} \mf u_\lambda = 0\; ,
\quad \text{for} \quad  k\ge 3
  \end{array}
\right.
\end{eqnarray*}
because $\bb J_+^* = - \bb J_-$ and because $\sigma$ has degree $2$.
The hierarchy starts at degree $2$ instead of $1$ because the degree
one equation is trivial. Indeed, by (\ref{eq:8}), $(\bb S + \bb L_1^c)
\pi_1 \mf u_\lambda =0$, $\bb J_- \pi_2 \mf u_\lambda =0$, so that the
degree one equation
$$
- \bb J_- \pi_{2} \mf u_\lambda + (\lambda-\bb S- \bb L_1^c) \pi_1
\mf u_\lambda = 0
$$
becomes $\pi_1 \mf u_\lambda = 0$. Hence $\pi_1 \mf u_\lambda$ plays
no role and we can set $\pi_1 \mf u_\lambda =0$.

Notice that we are using the same notation $\sigma$ for the local
function defined in (\ref{eq:11}) and its Fourier transform $\sigma
: \mc E\to\bb R$ which is given by
\begin{equation}
\label{eq:12}
\sigma (A) \;=\; \theta_j \Big\{
r_0 (e_j\cdot v) + \sum_{a=1}^d r_a (e_a\cdot v) (e_j\cdot v) \Big\}
\end{equation}
if $A = \{(0,v), (e_j,v)\}$ for some $v$ in $\mc V$, $1\le j\le d$,
and $\sigma (A) =0$, otherwise.

Consider the truncated resolvent equation up to the degree $n$:
\begin{equation}
\label{eq:5}
\left\{
\begin{array}{l}
\vphantom{\Big\{}
\bb J_+^* \pi_{3} \mf u_{\lambda,n} + (\lambda-\bb S - \bb L_1^c)
\pi_2 \mf u_{\lambda,n} =  \sigma \; , \\
\vphantom{\Big\{}
\bb J_+^* \pi_{k+1} \mf u_{\lambda,n} +
(\lambda- \bb S - \bb L_1^c) \pi_k \mf u_{\lambda,n}
- \bb J_+ \pi_{k-1} \mf u_{\lambda,n} =  0\; , \quad
3 \le k \le n-1\; , \\
\vphantom{\Big\{}
(\lambda- \bb S - \bb L_1^c) \pi_n \mf u_{\lambda,n}
- \bb J_+ \pi_{n-1} \mf u_{\lambda,n} = 0\; .
\end{array}
\right.
\end{equation}
We can solve  the final  equation  of \eqref{eq:5} by
$$
\pi_n \mf u_{\lambda,n} = (\lambda -\bb S - \bb L_1^c)^{-1}
\bb J_{+} \pi_{n-1} \mf u_{\lambda,n} \; .
$$
Substituting this into the equation of degree $n-1$, we have
$$
\pi_{n-1} \mf u_{\lambda,n} = \Big [ (\lambda -\bb S - \bb L_1^c)
+ \bb J_+^*  (\lambda -\bb S - \bb L_1^c) ^{-1} \bb J_{+} \Big]^{-1}
\bb J_+ \pi_{n-2} \mf u_{\lambda,n}\; .
$$
Solving iteratively we arrive at
$$
\pi_2 \mf u_{\lambda,n} \; =\; \mc T_n \sigma\; ,
$$
where the operators $\{\mc T_n, n\ge 2\}$ are defined inductively by
\begin{equation}
\label{eq:y1}
\mc T_2 = (\lambda - \bb S - \bb L_1^c)^{-1} \; ,\quad
\mc T_{n+1} \; =\; \Big\{ (\lambda - \bb S - \bb L_1^c) +
\bb J_{+}^* \mc T_n^{-1}  \bb J_{+} \Big \}^{-1}\; .
\end{equation}

The truncated equation represents the solution of $ (\lambda -\mc L_n)
\mf u_{\lambda,n} = \sigma$ and hence $\ll \pi_2 \mf u_{\lambda,n}
, \sigma\gg = \ll \sigma ,(\lambda-\mc L_n)^{-1} \sigma\gg$ so that
$$
\ll \sigma ,(\lambda-\mc L_n)^{-1} \sigma\gg \; =\; \ll \sigma, \mc T_n
\sigma \gg \; ,
$$
where, for example,
\begin{eqnarray*}
\!\!\!\!\!\!\!\!\!\!\!\!\!\!\!\!
&& \mc T_3 \; =\; \Big \{  (\lambda - \bb S - \bb L_1^c)  + \bb J_-
(\lambda - \bb S - \bb L_1^c)^{-1}  \bb J_{+} \Big \}^{-1} \; , \\
\!\!\!\!\!\!\!\!\!\!\!\!\!\!\!\!
&& \mc T_4 \; =\;
\Big [  (\lambda-\bb S - \bb L_1^c) +  \bb J_-
\Big \{ (\lambda-\bb S - \bb L_1^c) + \bb J_-
 (\lambda-\bb S - \bb L_1^c)^{-1}\bb J_{+} \Big \}^{-1}
\bb J_{+} \Big ]^{-1}\; .
\nonumber
\end{eqnarray*}

\begin{lemma}
\label{s3}
For each $\lambda>0$, $\ll \sigma , (\lambda - \mc L_{2k+1})^{-1} \sigma
\gg$ is an increasing sequence which converges to $\ll \sigma ,
(\lambda - \mc L^{ex} - \mc L_1^c)^{-1} \sigma \gg$ and $\ll \sigma ,
(\lambda - \mc L_{2k})^{-1} \sigma \gg$ is a decreasing sequence which
converges to $\ll \sigma , (\lambda - \mc L^{ex} - \mc L_1^c)^{-1}
\sigma \gg$.
\end{lemma}

\begin{proof}
Since $\lambda-\bb S - \bb L_1^c$ is positive, it is easy to show from
the definition of the sequence of operators $\mc T_n$ that $0\le \mc
T_3 \le \mc T_2$ and that $\mc T_m \le \mc T_n$ if $\mc T_{m-1} \ge
\mc T_{n-1}$. In particular, $\{\mc T_{2k} ,\, k\ge 1\}$ is a
decreasing sequence, $\{\mc T_{2k+1} ,\, k\ge 1\}$ is an increasing
sequence and $\mc T_{2k+1} \le \mc T_{2j}$ for any $k,j\ge 1$:
\begin{eqnarray}
\label{eq:6}
&& \ll \sigma, (\lambda-\mc L_3)^{-1} \sigma \gg \; \le\;
\ll \sigma, (\lambda- \mc L_5)^{-1} \sigma \gg \;\le\; \cdots  \\
&& \quad  \cdots \le \; \ll \sigma, (\lambda-\mc L_4)^{-1} \sigma \gg
\;\le\; \ll \sigma, (\lambda-\mc L_2)^{-1} \sigma\gg \; .
\nonumber
\end{eqnarray}

To check that $\ll \sigma, (\lambda-\mc L^{ex} - \mc L_1^c)^{-1} \sigma
\gg$ is in fact the limit of these upper and lower bounds we use the
variational formula.  For any matrix $M$, let $M_s$ denote the
symmetric part $(M+M^*)/2$.  The identity $ \big\{ [M^{-1}]_s
\big\}^{-1} =M^\ast (M_s)^{-1} M $ always holds, and thus we have
\begin{eqnarray}
\label{eq:7}
\!\!\!\!\!\!\!\!\!\!\!\!\! &&
\ll \sigma, (\lambda-\mc L^{ex} - \mc L_1^c)^{-1} \sigma \gg \\
\!\!\!\!\!\!\!\!\!\!\!\!\! && \quad
=\; \sup_{\mf f} \Big\{ 2\ll \sigma,\mf f\gg -
\ll  (\lambda-\mc L^{ex} - \mc L_1^c) \mf f,
(\lambda-\bb S - \bb L_1^c)^{-1}
(\lambda-\mc L^{ex} - \mc L_1^c)\mf f\gg \Big\}\; ,
\nonumber
\end{eqnarray}
where the supremum is carried over all finite supported functions $\mf
f:\mc E\to \bb R$.  Note that
\begin{eqnarray*}
\!\!\!\!\!\!\!\!\!\!\!\!\! &&
\ll (\lambda-\mc L^{ex} - \mc L_1^c)\mf f,
(\lambda-\bb S - \bb L_1^c)^{-1}
(\lambda-\mc L^{ex} - \mc L_1^c)\mf f\gg \\
\!\!\!\!\!\!\!\!\!\!\!\!\! && \qquad\qquad\qquad\qquad
=\; \ll \mf f, (\lambda-\bb S - \bb L_1^c) \mf f\gg +
\ll \bb A \mf f, (\lambda - \bb S - \bb L_1^c)^{-1} \bb A \mf f\gg\; ,
\end{eqnarray*}
where $\bb A = \bb J_+ + \bb J_-$.  Hence,
\begin{eqnarray*}
\!\!\!\!\!\!\!\!\!\!\!\!\! &&
\ll \sigma, (\lambda-\mc L^{ex} - \mc L_1^c)^{-1} \sigma\gg
=\; \sup_{\mf f} \inf_{\mf g}
\Big \{ 2\ll \sigma -\bb A^* \mf g, \mf f\gg \\
\!\!\!\!\!\!\!\!\!\!\!\!\! && \qquad\qquad\qquad\qquad
\;-\;
\ll \mf f, (\lambda-\bb S - \bb L_1^c) \mf
f\gg \;+\; \ll \mf g, (\lambda - \bb S - \bb L_1^c) \mf g\gg\Big\} \; .
\end{eqnarray*}
Let $a_n$ denote the supremum restricted to finite supported functions
in $\mf f$ in $L^2(\mc G_n)$, and $a^n$ denote the infimum restricted
to finite supported function $\mf g$ in $L^2(\mc G_n)$ so that
$a_n\uparrow\ll \sigma, (\lambda-\mc L^{ex} - \mc L_1^c)^{-1}
\sigma\gg$ and $a^n\downarrow\ll \sigma, (\lambda-\mc L^{ex} - \mc
L_1^c)^{-1} \sigma \gg$.  By straightforward computation one checks
that $a_n\le \ll \sigma , (\lambda-\mc L_{n+1})^{-1} \sigma\gg \le
a^n$, giving the desired result.
\end{proof}

In what follows we will present a general approach to the equations
(\ref{eq:5}) which, from (\ref{eq:6}) gives a nontrivial lower bound
on the diffusion coefficient.  Because it gives a sequence of upper
and lower bounds, the method has the potential to give the full
conjectured scaling of the diffusion coefficient.

\section{Removal of hard core}
\label{sec4}

From Lemma \ref{s3} of the previous section we have a lower bound at
degree three.  However, computations are complicated due to the hard
core exclusion.  We follow \cite{LQSY} to remove the hard core
restriction in the formulas and then perform explicit computations in
Fourier space. By removal of the hard core, we mean replacing
functions defined on $\mc E_n$ by symmetric functions defined on $E_n
= (\bb Z^d \times \mc V)^{n}$ and replacing operators acting on $\mc
E_n$ by operators acting on $E_n$.

We first identify a function $f: \mc E_n \to \bb R$ with a symmetric
function $f : E_n \to \bb R$. Denote by $\bs \omega_n = (\omega_1,
\dots, \omega_n)$, $\omega_i = (x_i,v_i)$, the points of $E_n$. For
$n\ge 1$, let
\begin{equation}
\label{eone}
E_{n,1} = \{ \bs\omega_n \, : \omega_i \not = \omega_j,
\text{ for } i \not = j \}
\end{equation}
and define
\begin{equation*}
f(\bs\omega_n) =
\begin{cases}
f (\{\omega_1,\dots,\omega_n\})
& {\rm if }\quad \bs\omega_n\in E_{n,1}\; , \\
0 & {\rm otherwise}\; .
\end{cases}
\end{equation*}

With the notation just introduced,
$$
E_{\mu_0}\Big [ \, \big ( \sum_{A\in\mc E_n} \mf f_A \Psi_A \big ) ^2 \Big ]
= \frac 1 {n! \, 4^n} \sum_{\bs \omega_n \in E_n} \mf f(\bs \omega_n)^2\; .
$$
For a function $f: E_n\to \bb R$, we shall use the same symbol
$\langle f \rangle $ to denote the expectation
$$
\frac 1 {n! \, 4^n} \sum_{\bs \omega_n \in E_n} f(\bs \omega_n)
$$
and write the inner product of two functions as $\langle f,g\rangle
= \langle fg \rangle$.  If $f$ and $g$ vanish on the complement of
$E_{n,1}$, this coincides with the inner product introduced before.
We also define, as before, $\ll f, g\gg = \sum_{x\in \bb Z^d} \langle
\tau_x f, g\rangle$.

Let $E = \cup_{n\ge 1} E_n$, $G_n = \cup_{1\le j\le n} E_j$.  We use
the same symbol $\pi_n$, $P_n = \sum_{1\le j\le n} \pi_j$ for the
projection onto $E_n$, $G_n$.  As before, there is a simple formula
for the inner product $\ll \, \cdot \, , \, \cdot \, \gg$. Consider
two finitely supported functions $f$, $g: E_n\to \bb R$. By
definition,
$$
n! \, 4^n \ll f , g \gg \; =\;
\sum_{\substack{\bs\omega_n\in E_n \\ z\in\bb Z^d}}
f(\bs \omega_n +z) g(\bs \omega_n) \; =\;
\sum_{\bs\omega_n\in E_n} \tilde f(\bs\omega_n) g(\bs\omega_n)\; ,
$$
where $\bs\omega_n + z = ((x_1+z, v_1), \dots, (x_n+z,v_n))$.

Denote by $\sim$ the equivalence relation on $E_n$ defined by
$\bs\omega_n \sim \bs\omega_n'$ if for all $1\le i\le n$, $v_i =v_i'$,
$x_i - x_i' =z$ for some $z$ in $\bb Z^d$.  Let $\tilde E_n = E_n
\big\vert_\sim$. Since summing over all sites in $E_n$ is the same as
summing over all equivalence classes and then over all elements of a
single class, the previous sum is equal to
$$
\sum_{\bs\omega_n\in \tilde E_n} \sum_{z\in \bb Z^d}
\tilde f(\bs\omega_n + z) g(\bs\omega_n+z) \; =\;
\sum_{\bs\omega_n\in \tilde E_n}
\tilde f(\bs\omega_n) \tilde g(\bs\omega_n)
$$
because $\tilde f(\bs\omega_n + z) = \tilde f(\bs\omega_n)$. It
remains to choose an element of each class.  This can be done by
fixing the last coordinate $x_n$ to be zero. In conclusion,
\begin{equation}
\label{eq:10}
n! \, 4^n \ll f , g \gg \; =\; \sum_{v\in\mc V}
\sum_{\bs\omega_{n-1}\in E_{n-1}}
f^*(\bs\omega_{n-1},v) g^*(\bs\omega_{n-1},v)\; ,
\end{equation}
where
\begin{equation}
\label{eq:9}
f^*(\bs\omega_{n-1},v) \; =\; \sum_{z\in \bb Z^d}
f((x_1+z, v_1), \dots, (x_{n-1}+z, v_{n-1}), (z,v))\; .
\end{equation}
Here again we see that the translations in the inner product
effectively reduce the degree of a function by one.

We derive now explicit formulas for the operators $\bb S$, $\bb A_+$
acting on symmetric functions of $E_n$. An elementary computation
shows that
\begin{equation*}
(\bb S f) (\bs \omega_n) \; =\; \gamma
\sum_{k = 1}^d \sum_{ i = 1}^ {n} \sum_{\iota = \pm }
\mb 1\{ \mb \nabla_{k,i}^\iota \bs \omega_n \in E_{n,1} \}
\nabla_{k,i}^\iota f (\bs \omega_n)
\end{equation*}
if $\bs \omega_n$ belongs to $E_{n,1}$ and $(\bb S f)(\bs \omega_n) =
0$ if $\bs \omega_n$ does not. Here, for $\iota = \pm$,
\begin{eqnarray*}
&& \mb \nabla_{k,i}^\iota \bs \omega_n \; =\;
(\omega_1, \dots , \omega_{i-1}, (x_{i}+ \iota e_k, v_i),
\omega_{i+1}, \dots, \omega_{n}) \; , \\
&&\quad
\nabla_{k,i}^\iota f (\bs \omega_n) \; =\;
f(\mb \nabla_{k,i}^\iota \bs \omega_n) - f(\bs \omega_n)\; .
\end{eqnarray*}
Note that $\bb S$ is the discrete Laplacian with Neumann boundary
condition on $E_{n,1}$. In the same way,
\begin{equation}
\label{Afdef}
(\bb J_+ f)(\bs \omega_n) \; =\;
\sum_{ i = 1}^ {n} \sum_{j \not =  i } \sum_{k=1}^d (e_k\cdot v_i)
\mb 1\{x_{j}+e_k = x_{i}, v_j=v_i \} \nabla_+^{i,j} f (\bs \omega_n)
\end{equation}
if $\bs \omega_n$ belongs to $E_{n,1}$ and $(\bb J_+ f)(\bs \omega_n)
= 0$ otherwise.  Here,
$$
\nabla_+^{i,j}f(\bs \omega_n) \; =\;  f(\bs \omega_n^i)
- f(\bs \omega_n^j)
$$
and the index $j$ in $\bs \omega_n^j$ indicates the absence of $
\omega_{j}$ in the vector $\bs \omega_n$: $\bs \omega_n^j = (\omega_1,
\dots, \omega_{j-1}, \omega_{j+1}, \dots, \omega_n)$. Finally, notice
that
\begin{eqnarray*}
\!\!\!\!\!\!\!\!\!\!\!\!\!\!\!\! &&
(\bb L_1^c f)(\bs \omega_n) \; =\; (1/4) \sum_{q\in \mc Q}
\sum_{j=1}^n i_q(v_j) \mb 1\{ x_k \not = x_j \text{ for }
k\not = j\} \\
\!\!\!\!\!\!\!\!\!\!\!\!\!\!\!\! &&  \qquad\qquad\qquad\qquad
\qquad\qquad
[ f(\sigma_{j,v'} \bs \omega_n) +
f(\sigma_{j,w'} \bs \omega_n) - f(\sigma_{j,v} \bs \omega_n) -
f(\sigma_{j,w} \bs \omega_n)]\;,
\end{eqnarray*}
if $\bs \omega_n$ belongs to $E_{n,1}$ and $(\bb L_1^c f)(\bs
\omega_n) = 0$ otherwise.  Here, $\sigma_{j,u} \bs \omega_n =
(\omega_1, \dots, \omega_{j-1}$, $(x_j,u), \omega_{j+1}, \dots,
\omega_n)$.

We now extend the operators $\bb S$, $\bb J_+$ to symmetric functions
not necessarily vanishing on $E_{n,1}$ by formulas analogous to the
ones above, except that we drop some indicator functions. Let $S$,
$J_+$ and $L_1^c$ be the operators defined by:
\begin{eqnarray*}
\!\!\!\!\!\!\!\!\!\!\!\!\!\!\!\! &&
(S F) (\bs \omega_{n}) \;=\; \gamma (\Delta F) (\bs \omega_{n}) \;=\;
\gamma \sum_{ i = 1}^ {n} \sum_{\iota = \pm } \sum_{k = 1}^d
(\nabla_{k,i}^\iota F) (\bs \omega_{n})\; , \\
\!\!\!\!\!\!\!\!\!\!\!\!\!\!\!\! &&  \quad
(J_+ F) (\bs \omega_{n}) \;=\;  \sum_{ i = 1}^ {n}
\sum_{j \not =  i } \sum_{k=1}^d (e_k\cdot v_i)
\mb 1\{ x_{j} +e_k= x_{i}, v_i=v_j \}
{\nabla}^{i,j}_+ F(\bs \omega_{n}) \quad\text{and} \\
\!\!\!\!\!\!\!\!\!\!\!\!\!\!\!\! && \qquad
(L_1^c f)(\bs \omega_n) \; =\; (1/4) \sum_{q\in \mc Q}
\sum_{j=1}^n i_q(v_j) \mb 1\{ x_k \not = x_j \text{ for }
k\not = j\} \\
\!\!\!\!\!\!\!\!\!\!\!\!\!\!\!\! &&  \qquad\qquad\qquad\qquad
\qquad
[ f(\sigma_{j,v'} \bs \omega_n) +
f(\sigma_{j,w'} \bs \omega_n) - f(\sigma_{j,v} \bs \omega_n) -
f(\sigma_{j,w} \bs \omega_n)]\;.
\end{eqnarray*}

Notice that $\langle L_1^c F \rangle = \langle J_+ F \rangle = 0$ if
$\langle |F|\rangle < \infty$ and hence the counting measure is
invariant. Let
\begin{equation*}
L= S + L_1^c + J_+\; .
\end{equation*}
and denote by $L_n = P_n L P_n$ the restriction of $L$ to $G_n$.
Following section 4 in \cite{LQSY}, we prove the next result which
permits to avoid the hard core interaction of the exclusion.

\begin{proposition}
\label{s4}
In dimension $d=2$, there exists a finite constant $C_0$ such that
\begin{equation*}
\frac{1}{C_0 n^{6}} {\ll} \sigma , (\lambda - L_n)^{-1} \sigma {\gg}
\;\le\; {\ll} \sigma,(\lambda -  {\mathcal L}_n)^{-1} \sigma {\gg}
\;\le\; C_0 n^4 {\ll} \sigma, (\lambda -L_n)^{-1} \sigma {\gg} \; .
\end{equation*}
for all $\lambda>0$.
\end{proposition}

The proof of this proposition is similar to the one of Lemma 3.1 in
\cite{LQSY} and therefore omitted. The main difference is to prove
Lemma 4.3 in \cite{LQSY} with $S + L_1^c$, $\bb S + \bb L_1^c$ in
place of $S$, $\bb S$ and this is elementary.

The special case $n=3$ combined with Lemmas \ref{s1}, \ref{s3} tells
us that
\begin{equation*}
\ll \sigma, (\lambda -{\mathcal L})^{-1} \sigma \gg \;\ge\;
C_0 \ll \sigma, (\lambda - L_3)^{-1} \sigma  \gg \; .
\end{equation*}

\section{Fourier computations}

To bound below $\ll \sigma, (\lambda - L_3)^{-1} \sigma \gg$, define
the Fourier transform of a function $\mf f:(\bb Z^{d} \times \mc
V)^n \to \bb R$ by
$$
\widehat{\mf f} \, (\bs p_{n}, \bs v_n) \;=\; \sum_{\bs x_n\in\bb Z^{nd}}
e^{-i \bs x_n\cdot \bs p_n} \mf f (\bs x_n, \bs v_n)
$$
for ${\bs p}_n \in (\bb R^d/2\pi \bb Z^d)^n$. Here we represented $\bs
w_n =(w_1, \dots , w_n)$, $ w_i=(x_i,v_i)$, as $(\bs x_n, \bs v_n)$
with $\bs x_n = (x_1, \dots, x_n) \in (\bb Z^d)^n$, $\bs v_n = (v_1,
\dots, v_n) \in \mc V^n$.

An elementary computation together with (\ref{eq:10}) shows that for
any local functions $f$, $g$ of degree $n$,
$$
\ll f, g \gg \;=\; \frac 1{(2\pi)^{(n-1)d} n! \, 4^n}
\sum_{\bs v_n \in \mc V^n} \int_{[-\pi, \pi]^{(n-1)d}}
\widehat{\mf f^*} (\bs p_{n-1}, \bs v_n)
\widehat{\mf g^*} (\bs p_{n-1}, \bs v_n)\, d\bs p_{n-1}\;.
$$
In this formula, $\mf f = \bb T f$, $\mf f^*$ is defined by
(\ref{eq:9}) and $\widehat{\mf f^*}$ is the Fourier transform of $(\mf
f^*) (\bs x_{n-1}, \bs v_n)$. Expressing $\widehat{\mf f^*}$ in terms
of $\widehat{\mf f}$ we further obtain that
$$
\ll f, g \gg \;=\; \frac 1{(2\pi)^{(n-1)d} n! \, 4^n}
\sum_{\bs v_n \in \mc V^n} \int_{\substack{[-\pi,
    \pi]^{nd} \\ \sum_{1\le j\le n} p_j =0}}
\widehat{\mf f} (\bs p_n, \bs v_n) \,
\widehat{\mf g} (\bs p_{n}, \bs v_n)\, d\bs p_{n}\;.
$$

Fix a symmetric function $\mf f:(\bb Z^{d} \times \mc V)^n \to \bb
R$.  The Fourier transform of the discrete Laplacian acting on $\mf f$
is given by
$$
- \widehat {\Delta \mf f} (\bs p_{n}, \bs v_n) \;=\;
\widehat {\mf f} (\bs p_{n}, \bs v_n) W(\bs p_{n})\;,
$$
where
$$
W(\bs p_{n}) \;=\; \sum_{j=1}^n W(p_j)
\;=\; \sum_{j=1}^n \sum_{k=1}^d
\{ 1 - \cos (p_j \cdot e_k) \}
$$
if $\bs p_{n} = (p_1, \dots p_n)$.  Notice that we are using the
same notation $W(\cdot)$ for slightly different objects. Moreover, for
$n=2$, a straightforward computation shows that
\begin{eqnarray*}
\!\!\!\!\!\!\!\!\!\!\!\!\!\!\! &&
\widehat { J_+ \mf f} \, ( p _1,  p _2,  p _3, \bs v_3) \;=\;
- i \sum_{j=1}^d \sum_{\sigma}  \mb 1\{ v_{\sigma_1} =
v_{\sigma_2}\} \, (e_j \cdot v_{\sigma_1}) \;\times \\
\!\!\!\!\!\!\!\!\!\!\!\!\!\!\! && \qquad\qquad \qquad \qquad
\{ \sin (e_j \cdot p _{\sigma_1}) + \sin (e_j \cdot
p _{\sigma_2}) \} \, \widehat {\mf f} ( p _{\sigma_1}+ p _{\sigma_2},  p
_{\sigma_3}; v_{\sigma_1} , v_{\sigma_3})\; ,
\end{eqnarray*}
where $\sigma$ runs over all permutations of degree three.

We are now in a position to state the first estimate based on Fourier
arguments.

\begin{lemma}
\label{s5}
Fix a symmetric function $\mf f:(\bb Z^d \times \mc V)^2 \to \bb R$.
There exists a finite constant $C_0$ such that in dimension $2$,
\begin{eqnarray*}
\!\!\!\!\!\!\!\!\!\!\!\!\!\!\! &&
\ll  J_+ \mf f,  (\lambda - \Delta_3)^{-1} J_+ \mf f \gg \;\le\; \\
\!\!\!\!\!\!\!\!\!\!\!\!\!\!\! && \quad
C_0 \, \sum_{v_1, v_2\in\mc V}
\int_{[-\pi, \pi]^2} dp\,  W(p)
\, |\log (\lambda +  W(p) )| \,
|\, \widehat {\mf f} (p,-p; v_1,v_2)\, |^2 \; .
\end{eqnarray*}
\end{lemma}

The proof of this result is similar to the one of Lemma 3.2 in
\cite{LQSY} and therefore omitted.

For $0\le j\le 2$, let $\bb I_j = \bb T I_j$, $0\le j\le 2$ be the
symmetric functions in $E_n$ associated to the conserved quantities.
An elementary computation shows that
$$
\bb I_0(v) \;=\; 1 \; , \quad \bb I_1(v) \;=\; e_1 \cdot v
\; , \quad \bb I_2(v) \;=\; e_2 \cdot v\;.
$$

Let $Q$ be the non positive symmetric matrix corresponding to the
operator $L^c_1$ acting on functions depending only on one site.
Notice that $\bb I_{a}$, $a=0$, $1$, $2$, are eigenvectors of $Q$.
Denote by $\{\bb I_{a} ,\, 3\leq a\leq |\mc V| -1\}$ the other
eigenvectors of $Q$ and by $\{q_{a} ,\, 0\leq a\leq |\mc V|-1\}$ the
corresponding eigenvalues. Since $\bb I_0$, $\bb I_1$, $\bb I_2$ are
associated to conserved quantities, $q_0 = q_1 = q_2=0$. Let $\mc O$
be the orthogonal matrix which diagonalizes $Q$:
\begin{equation*}
\mc D = \mc O^* Q \mc O\;.
\end{equation*}
For $n\ge 1$, let
\begin{equation*}
Q_n = \sum_{i=1}^n I \otimes \cdots \otimes Q \otimes \cdots \otimes
I\;.
\end{equation*}
Since $L^c_1$ has an indicator function, we have that $0 \leq - L^c
\leq - Q_n$. Moreover, $Q_n$ can be diagonalized by $\mc O_n = \mc
O^{\otimes n}$: $\mc D_n = \mc O_n^* Q_n \mc O_n$, where $\mc D_n =
\sum I \otimes \cdots \otimes \mc D \otimes \cdots \otimes I$.  Notice
that the Laplacian commutes with these matrices.

As in section \ref{sec3}, we can represent $(\lambda - L_3)^{-1}
\sigma$ in terms of the operators $S$, $J_+$ and $L_1^c$ to obtain that
$$
\ll \sigma, (\lambda - L_3)^{-1} \sigma \gg
\;=\; \ll \sigma , \big\{ \lambda - \Delta - L^c_1 + J_+^*(\lambda -
\Delta -  L^c_1)^{-1} J_+\big\}^{-1} \sigma \gg\; .
$$
Since $0 \le -L_c^1 \le - Q_n$, the previous scalar product is
bounded below by
$$
\ll \sigma,  \big\{ \lambda - \Delta_2 - Q_2 +
J_+^*(\lambda - \Delta_3)^{-1} J_+\big\}^{-1} \sigma\gg\;.
$$
Recalling that $Q_2 = \mc O_2 \mc D_2 \mc O_2^*$ and that $\Delta$
commutes with the operators $\mc O_2$, we rewrite the previous
expression as
\begin{equation*}
\ll \mc O_2^*  \sigma,  \big\{ \lambda - \Delta_2 - \mc D_2 + \mc O_2^*
J_+^*(\lambda - \Delta_3)^{-1} J_+ \mc O_2 \big\}^{-1}
\mc O_2^*\sigma\gg\;.
\end{equation*}
To keep notation simple, let $\Omega = \{ \lambda - \Delta_2 - \mc D_2
+ \mc O_2^* J_+^*(\lambda - \Delta_3)^{-1} J_+ \mc O_2 \}^{-1}$ and
denote by $\bar{\pi}$ the projection onto the eigenspace corresponding
to the zero eigenvalue of $\mc D_2$. By Schwarz inequality, for any
function $H$,
$$
\ll H , \Omega H \gg \;\ge\; (1/2) \ll \bar{\pi} H , \Omega \bar{\pi}
H \gg \; -\; \ll (1-\bar{\pi}) H , \Omega (1-\bar{\pi}) H \gg\;.
$$

We claim that $\ll (1-\bar{\pi}) \mc O_2^* \sigma , \Omega
(1-\bar{\pi}) \mc O_2^* \sigma \gg$ is bounded by a finite constant.
Indeed, by definition of $\Omega$,
$$
\ll (1-\bar{\pi}) \mc O_2^* \sigma , \Omega (1-\bar{\pi})
\mc O_2^* \sigma \gg \;\le\;
\ll \mc (1-\bar{\pi}) O_2^*  \sigma,  (-\mc D_2)^{-1} (1-\bar{\pi})
\mc O_2^*\sigma\gg\;.
$$
Since $(1-\bar{\pi})$ is the projection on the positive eigenvalues
of $\mc D_2$, the previous expression is less than or equal to a
finite constant which depends on a lower bound for the positive
eigenvalues of $Q$ and an upper bound for $\sigma$.

The following lemma is needed to estimate $\ll \bar{\pi} \mc O_2^*
\sigma , \Omega \bar{\pi} \mc O_2^* \sigma \gg$. Recall that $\sigma =
\sum_{0\le a\le d} \sum_{1\le j\le d} r_a \theta_j \sigma^a_j$.

\begin{lemma}
\label{s6}
$\bar{\pi}\mc O^*_2 \sigma = 0$ if and only if $\theta =0$ or $r=0$.
\end{lemma}

\begin{proof}
Assume that $\bar{\pi} \mc O^*_2 \sigma =0$ and assume without loss of
generality that $\theta_1 \neq 0$. Fix the component $\{0,e_1\}$.
Since $\bar{\pi}$ is the projection of the eigenspace of $\mc D_2$
associated to the zero eigenvalues, $\bar{\pi} \mc O^*_2 \sigma =0$ if
and only if the scalar product of $\sigma$ with $(\bb I_a, \bb I_b)$
vanishes for $0\le a,b\le 2$.

An elementary computation based on the explicit expression (\ref{eq:12})
for $\sigma$ shows that
$$
\sigma \cdot (\bb I_a, \bb I_b) \;=\; \theta_1
\sum_{v\in\mc V} \bb I_a (v) \, \bb I_b (v) \big\{ r_0 (e_1 \cdot v) +
r_1 (e_1 \cdot v)^2 + r_2 (e_1 \cdot v) (e_2 \cdot v) \big\}\;.
$$
Setting $a = 0$ and $b=1$ we get that
\begin{equation*}
\theta_1 \sum_{v\in \mc V} \big\{ r_0 (e_1\cdot v)^2  +
r_1 (e_1\cdot v)^3  + r_2 (e_1\cdot v)^2 (e_2\cdot v)
\big\}\; .
\end{equation*}
The last two terms vanish due to the symmetry of $\mc V$. Hence $r_0 =
0$. Repeating the same argument with $a = b = 1$ (resp. $a=1, b=2$),
we obtain that $r_1=0$ (resp. $r_2=0$). This concludes the proof of
the lemma.
\end{proof}

We are now in a position to bound below
$$
\ll \bar{\pi} \mc O_2^* \sigma , \{ \lambda - \Delta_2 - \mc D_2
+ \mc O_2^* J_+^*(\lambda - \Delta_3)^{-1} J_+ \mc O_2 \}^{-1}
\bar{\pi} \mc O_2^* \sigma \gg\;.
$$

Fix a function $G$ in the image of $\bar \pi$. By the variational
formula for the $H_{-1}$ norm, the previous scalar product with $G$ in
place of $\bar{\pi} \mc O_2^* \sigma$ is equal to
$$
\sup_{F} \Big\{ 2\ll G , F\gg
- \ll F, [\lambda - \Delta_2 -\mc D_2] F \gg - \ll \mc O_2 F,
J_{+}^* (\lambda-\Delta_3)^{-1}  J_{+} \mc O_2 F\gg \Big\}\; ,
$$
where the supremum is performed over all local functions $F$.
Restricting this supremum to functions in the image of $\bar{\pi}$ we
obtain a lower bound. In this case we may remove the operator $\mc
D_2$ because $\bar{\pi}$ is the projection onto the space associated
to the zero eigenvalues of $\mc D_2$.

Lemma \ref{s5} and a straightforward computation show that
$$
\ll \mc O_2 F, J_{+}^* (\lambda-\Delta_3)^{-1}  J_{+} \mc O_2 F\gg
$$
is bounded above by the right hand side of the statement of Lemma
\ref{s5}. Now, repeating the arguments presented in the proof of Lemma
3.3 in \cite{LQSY}, we obtain that the previous variational formula is
bounded below by
\begin{equation*}
C_0 \sum_{\bs v_2} \int_{[-\pi,\pi)^2}
\frac { | \widehat{G}(p,-p, \bs v_2)|^2}
{\lambda + 2 W(p)  + C_1 W(p) | \log\{\lambda + W(p)\}| }
\, dp
\end{equation*}
for some finite constants $C_0$, $C_1$. We used here that $G$ is in the
image of $\bar{\pi}$ because we needed the function $F$ which
maximizes the supremum to be in the image of $\bar{\pi}$. The factor
$2$ appeared because $W(p,-p) = 2 W(p)$.

Replace $G$ by $\bar{\pi} \mc O_2^* \sigma$. It follows from the
previous lemma and elementary computations that $\sum_{\bs v_2} |
\widehat{\bar{\pi} \mc O_2^* \sigma}(p,-p, \bs v_2)|^2$ is bounded
below by a positive constant close to the origin.  The previous
expression is therefore bounded below by

\begin{equation*}
C_0 \int_{[-\epsilon,\epsilon)^2}
\frac {1}
{\lambda + 2 W(p)  + C_1 W(p) | \log\{\lambda + W(p)\}| }
\, dp
\end{equation*}
for some $\epsilon>0$.  Changing to polar coordinates, since there is
no angle dependence, we get an integral of the form
\begin{equation*}
\int_0^{\rho} \frac{rdr}{\lambda + r^2(2 - C_1 \log(\lambda+ r^2))} \sim
\int_0^{\rho} \frac{-r \, dr}{(\lambda + r^2) \log(\lambda+ r^2)}
\end{equation*}
for some $\rho>0$.  By the change of variables $u = -\log(\lambda
+ r^2)$ we finally get
\begin{equation*}
\sim \int_{-\log{\rho^2}}^{-\log\lambda} \frac{du}{u} \sim
\log{\log\lambda^{-1}} \;.
\end{equation*}
This proves Lemma \ref{s-1}.

\bigskip

Recall the recursive relation \eqref{eq:y1}. If we denote the limit
of $\mc T_n = \mc T$, then it satisfies the equation
\begin{equation*}
\mc T \; =\;
\Big [  (\lambda-\bb S - \bb L_1^c) +  \bb J_-
\mc T ^{-1} \bb J_{+} \Big ]^{-1}\;.
\end{equation*}
We now assume that the dispersion relation of ${\mc T}^{-1}$ is given by
$$
\Big \{ \sum_{j=1}^\infty W(p_j) \Big \}
\, \Big |\log \Big (\lambda +  \Big \{ \sum_{j=1}^\infty W(p_j)
\Big \} \Big) \, \Big |^\kappa
$$
for some $\kappa>0$.  Write $\sum_{j=1}^\infty W(p_j) = u + W(p_1)$
where $u=\sum_{j=2}^\infty W(p_j)$. Suppose that we are interested in
$u\ge \lambda$.  Then $\bb J_- \mc T ^{-1} \bb J_{+}$ is approximately
given by
\begin{equation}
\label{3}
\int_{[-\epsilon,\epsilon)^2}
\frac {1}
{ u + W(p)  +  (u+ W(p))| \log\{u + W(p)\}|^\kappa }
\, dp\sim u |\log u|^{1-\kappa}
\end{equation}
For $u \ge \lambda$ we have
$$\mc T ^{-1} = (\lambda-\bb S - \bb L_1^c) +  \bb J_-
\mc T ^{-1} \bb J_{+} \sim   \bb J_-
\mc T ^{-1} \bb J_{+}
$$
Thus we have $\kappa= 1- \kappa$ and this gives the value
$\kappa=1/2$.  This is precisely the exponent derived in
\cite{AWG,FNS}.

Notice that the exponent is different from the ASEP which takes the
value $\kappa=2/3$ \cite{Y}. This is because that the dispersion law
is changed only in one direction for ASEP. In the fluid model
considered here, the collision operator spread the dispersion law to
all direction. This causes the exponent on the right side of \eqref{3}
to be $1-\kappa$, compared with $1-\kappa/2$ in \cite{Y}. This sketch
is certainly not rigorous since it is not even clear where the
operator $\mc T$ should be defined. Though one can try to prove upper
and lower bound as in \cite{Y}, it is not clear that the off-diagonal
terms can be controlled. However, the exponent $1/2$ seems to be very
convincing.
\medskip

\noindent{\bf Acknowledgments.} C. Landim was partially supported
by CNPq grant 474626/03-2 and Projeto PRONEX ``Probabilidade e
Processos Estoc\'asticos''. Part of this work was done while J.
Ramirez was a postdoctor at the Cornell University and also during
visits to the New York University; he would like to thank the
hospitality of these institutions. H.-T. Yau was partially
supported by the NSF grant DMS-0307295 and a MacArthur Fellowship.


\begin{thebibliography}{60}
  
\bibitem{AW} B. J. Alder, T.E. Wainwright: {Decay of the Velocity
    Autocorrelation Function.} Phys. Rev. A {\bf 1}, 18-21 (1970)

\bibitem{AWG} B. J. Alder, T.E. Wainwright, D. Gass: {Decay of Time
    Correlations in Two Dimensions.} Phys. Rev. A {\bf 4}, 233-237
  (1971)

\bibitem{BKS} H. van Beijeren, R. Kutner, H. Spohn: { Excess noise for
    driven diffusive systems.} Phys. Rev. Lett. {\bf 54}, 2026-2029
  (1985)

\bibitem{EMYa} Esposito, R.; Marra, R.; Yau, H. T: Navier-Stokes
  Equations for Stochastic Particle Systems on the Lattice.  Comm.
  Math. Phys. {\bf 182}, 395--456, (1996).
  
\bibitem{FNS} D. Forster, D. Nelson, M. Sthephen: {Large-distance and
    Long-time Properties of a Random Stirred Fluid.} Phys. Rev. A {\bf
    16}, 732-749 (1970)

\bibitem{HF} M. van der Hoef, D. Frenkel: {Evidence for
    Faster-than-$t^{-1}$ Decay of the Velocity Autocorrelation
    function in a 2D Fluid.}  Phys. Rev. Lett. {\bf 66}, 1591-1594
  (1991)

\bibitem{J} K. Johansson: {Shape fluctuations and random matrices}.
  Comm.  Math. Phys. {\bf 209}, 437-476 (2000)

\bibitem{kl} C. Kipnis, C. Landim; {\it Scaling Limit of Interacting
    Particle Systems}, Grundlheren der mathematischen Wissenschaften
  {\bf 320}, Springer-Verlag, Berlin, New York, (1999).
  
\bibitem{lov6} C. Landim, S. Olla, S. R. S. Varadhan; On viscosity and
  fluctuation--dissipation in exclusion processes.  J. Stat. Phys.
  {\bf 115}, 323-363, (2004) \smallskip

\bibitem{loy} C. Landim, S. Olla, H. T. Yau; Some properties of the
  diffusion coefficient for asymmetric simple exclusion processes.
  Annals Probab. {\bf 24}, 1779--1807, (1996).

\bibitem{LQSY} Landim, C.; Quastel, J.; Salmhofer, M.;Yau, H. T.
  Superdiffusivity of asymmetric exclusion process in dimensions one
  and two. Comm. Math. Phys. {\bf 244}, 455-481 (2003)

\bibitem{LY} Landim, C.; Yau, H. T. Fluctuation-dissipation equation
  of asymmetric simple exclusion processes. Probab. Th.  Rel. Fields
  {\bf 108}, 321--356, (1997).
  
\bibitem{QY} Quastel, J.; Yau, H.-T. Lattice gases, large deviations,
  and the incompressible Navier-Stokes equations. Ann. of Math. {\bf
    148}, 51--108, (1998).

\bibitem{SVY} Sethuraman, S.; Varadhan, S. R. S.; Yau, H.T. :
  Diffusive limit of a tagged particle in asymmetric simple exclusion
  processes. Comm. Pure Appl. Math. 53 (2000), no. 8, 972--1006.

  
\bibitem{Y} Yau, H. T.: $(\log t)^{2/3}$ law of the two dimensional
  asymmetric simple exclusion process, to appear in Annals Math.

\end{thebibliography}
\end{document}